\documentclass[reqno,11pt]{amsart}%
\usepackage{amsmath}
\usepackage{graphicx}
\usepackage{amsxtra}%
\usepackage{amsfonts}%
\usepackage{amssymb}
\theoremstyle{plain}
\newtheorem{theorem}{Theorem}[section]
\newtheorem{lemma}[theorem]{Lemma}
\newtheorem{proposition}[theorem]{Proposition}
\newtheorem{corollary}[theorem]{Corollary}

\theoremstyle{definition}

\newtheorem{example}[theorem]{Example}

\theoremstyle{remark}
\newtheorem{remark}[theorem]{Remark}

\numberwithin{equation}{section}
\begin{document}
\title{Truncated $K$-moment problems \\in several variables}
\author{Ra\'{u}l E. Curto}
\address{Department of Mathematics, The University of Iowa, Iowa City, Iowa 52242}
\email{rcurto@math.uiowa.edu}
\urladdr{http://www.math.uiowa.edu/\symbol{126}rcurto/}
\author{Lawrence A. Fialkow}
\address{Department of Computer Science, State University of New York, New Paltz, NY 12561}
\email{fialkowl@newpaltz.edu}
\urladdr{http://www.mcs.newpaltz.edu/faculty/fialkow.html}
\thanks{The first-named author's research was partially supported by NSF Research
Grant DMS-0099357. \ The second-named author's research was partially
supported by NSF Research Grant DMS-0201430. \ The second-named author was
also partially supported by the State University of New York at New Paltz
Research and Creative Projects Award Program.}

\begin{abstract}
Let $\beta\equiv\beta^{(2n)}$ be an $N$-dimensional real multi-sequence of
degree $2n$, with associated moment matrix $\mathcal{M}(n)\equiv
\mathcal{M}(n)(\beta)$, and let $r:=\operatorname*{rank}\mathcal{M}(n)$. \ We
prove that if $\mathcal{M}(n)$ is positive semidefinite and admits a
rank-preserving moment matrix extension $\mathcal{M}(n+1)$, then
$\mathcal{M}(n+1)$ has a unique representing measure $\mu$, which is
$r$-atomic, with $\operatorname*{supp}\mu$ equal to $\mathcal{V}%
(\mathcal{M}(n+1))$, the algebraic variety of $\mathcal{M}(n+1)$. \ Further,
$\beta$ has an $r$-atomic (minimal) representing measure supported in a
semi-algebraic set $K_{\mathcal{Q}}$ subordinate to a family $\mathcal{Q}%
\equiv\{q_{i}\}_{i=1}^{m}\subseteq\mathbb{R}[t_{1},...,t_{N}]$ if and only if
$\mathcal{M}(n)$ is positive semidefinite and admits a rank-preserving
extension $\mathcal{M}(n+1)$ for which the associated localizing matrices
$\mathcal{M}_{q_{i}}(n+[\frac{1+\deg q_{i}}{2}])$ are positive semidefinite
$(1\leq i\leq m)$; in this case, $\mu$ (as above) satisfies
$\operatorname*{supp}\mu\subseteq K_{\mathcal{Q}}$, and $\mu$ has precisely
$\operatorname*{rank}\mathcal{M}(n)-\operatorname*{rank}\mathcal{M}_{q_{i}%
}(n+[\frac{1+\deg q_{i}}{2}])$ atoms in $\mathcal{Z}(q_{i})\equiv\left\{
t\in\mathbb{R}^{N}:q_{i}(t)=0\right\}  $, $1\leq i\leq m$. \ 

\end{abstract}
\maketitle

\section{\label{Int}Introduction}

Given a finite real multisequence $\beta\equiv\beta^{(2n)}=\left\{  \beta
_{i}\right\}  _{i\in\mathbb{Z}_{+}^{N},\;\left|  i\right|  \leq2n}$, and a
closed set $K\subseteq\mathbb{R}^{N}$, the \textit{truncated }$K$%
\textit{-moment problem} for $\beta$ entails determining whether there exists
a positive Borel measure $\mu$ on $\mathbb{R}^{N}$ such that
\begin{equation}
\beta_{i}=\int_{\mathbb{R}^{N}}t^{i}\,d\mu(t),\qquad i\in\mathbb{Z}_{+}%
^{N},\;\left|  i\right|  \leq2n, \label{IntroEq1}%
\end{equation}
and%
\begin{equation}
\operatorname*{supp}\mu\subseteq K; \label{IntroEq2}%
\end{equation}
a measure $\mu$ satisfying (\ref{IntroEq1}) is a \textit{representing measure
}for $\beta$; $\mu$ is a $K$\textit{-representing measure} if it satisfies
(\ref{IntroEq1}) and (\ref{IntroEq2}).

In the sequel, we characterize the existence of a finitely atomic
$K$-represen\-ting measure having the fewest possible atoms, in the case when
$K$ is semi-algebraic. This is the case where $\mathcal{Q}\equiv\left\{
q_{i}\right\}  _{i=1}^{m}\subseteq\mathbb{R}^{N}\left[  t]\equiv
\mathbb{R[}t_{1},\ldots,t_{N}\right]  $ and $K=K_{\mathcal{Q}}:=\left\{
(t_{1},\ldots,t_{N})\in\mathbb{R}^{N}:q_{i}(t_{1},\ldots,t_{N})\geq0,1\leq
i\leq m\right\}  $. Our existence condition (Theorem \ref{IntroTheorem1}
below) is expressed in terms of positivity and extension properties of
the\textit{ moment matrix} $\mathcal{M}(n)\equiv\mathcal{M}^{N}(n)(\beta)$
associated to $\beta$, and in terms of positivity of the \textit{localizing
matrix} $\mathcal{M}_{q_{i}}$ corresponding to each $q_{i}$ (see below for
terminology and notation). \ In Theorem \ref{thm12} we provide a procedure for
computing the atoms and densities of a minimal representing measure in any
truncated moment problem (independent of $K$).

If $\mu$ is a representing measure for $\beta$ (or, as we often say, a
representing measure for $\mathcal{M}(n)$), then $\operatorname*{card}%
\operatorname*{supp}\mu\geq$ $\operatorname*{rank}\mathcal{M}(n)$; moreover,
there exists a $\operatorname*{rank}\mathcal{M}(n)$-atomic (minimal)
representing measure for $\beta$ if and only if $\mathcal{M}(n)$ is positive
semidefinite ($\mathcal{M}(n)\geq0$) and $\mathcal{M}(n)$ admits a
rank-preserving (or \textit{flat}) extension to a moment matrix $\mathcal{M}%
(n+1)$; in this case, $\mathcal{M}(n+1)$ admits unique successive flat moment
matrix extensions $\mathcal{M}(n+2),$ $\mathcal{M}(n+3),\ldots$ (Theorem
\ref{newthm}). \ For $1\leq i\leq m$, suppose $\deg q_{i}=2k_{i}$ or
$2k_{i}-1$; relative to $\mathcal{M}(n+k_{i})$ we have the \textit{localizing
matrix} $\mathcal{M}_{q_{i}}(n+k_{i})$ (cf. Section \ref{Loc}).

Our two main results, which follow, characterize the existence of
$\operatorname*{rank}M(n)$-atomic (minimal) $K_{\mathcal{Q}}$-representing
measures for $\beta$ and show how to compute the atoms and densities of such measures.

\begin{theorem}
\label{IntroTheorem1} An $N$-dimensional real sequence $\beta\equiv
\beta^{(2n)}$ has a $\operatorname*{rank}\mathcal{M}(n)$-atomic representing
measure supported in $K_{\mathcal{Q}}$ $\ $if and only if $\mathcal{M}%
(n)\geq0$ and $\mathcal{M}(n)$ admits a flat extension $\mathcal{M}(n+1)$ such
that $\mathcal{M}_{q_{i}}(n+k_{i})\geq0$ $(1\leq i\leq m)$. In this case,
$\mathcal{M}(n+1)$ admits a unique representing measure $\mu$, which is a
$\operatorname*{rank}\mathcal{M}(n)$-atomic (minimal) $K_{\mathcal{Q}}%
$-representing measure for $\beta$; moreover, $\mu$ has precisely
$\operatorname*{rank}\mathcal{M}(n)-\operatorname*{rank}\mathcal{M}_{q_{i}%
}(n+k_{i})$ atoms in $\mathcal{Z}(q_{i})\equiv\left\{  t\in\mathbb{R}%
^{N}:q_{i}(t)=0\right\}  $, $1\leq i\leq m$.
\end{theorem}%

\noindent
The uniqueness statement in Theorem \ref{IntroTheorem1} actually depends on
our next result, which provides a concrete procedure for computing the measure
$\mu$. \ As described in Section \ref{Mom}, the rows and columns of
$\mathcal{M}(n)$ are indexed by the lexicographic ordering of the monomials
$t^{i}\;(i\in\mathbb{Z}_{+}^{N},\;\left|  i\right|  \leq n)$, and are denoted
by $T^{i}\;(\left|  i\right|  \leq n)$; a dependence relation in the columns
of $\mathcal{M}(n)$ may thus be expressed as $p(T)=0$ for a suitable
$p\in\mathbb{R}^{N}[t]$ with $\deg p\leq n$. \ We define the \textit{variety}
of $\mathcal{M}(n)$ by $\mathcal{V}(\mathcal{M}(n)):=\bigcap_{\substack{p\in
\mathbb{R}^{N}[t],\deg p\leq n\\p(T)=0}}\mathcal{Z}(p)$, where $\mathcal{Z}%
(p):=\{t\in\mathbb{R}^{N}:p(t)=0\}$. \ Let $r:=\operatorname*{rank}%
\mathcal{M}(n)$ and let $\mathcal{B}\equiv\{T^{i_{k}}\}_{k=1}^{r}$ denote a
maximal linearly independent set of columns of $\mathcal{M}(n)$. \ For
$\mathcal{V}\equiv\{v_{j}\}_{j=1}^{r}\subseteq\mathbb{R}^{N}$, let
$W_{\mathcal{B},\mathcal{V}}$ denote the $r\times r$ matrix whose entry in row
$k$, column $j$ is $v_{j}^{i_{k}}\;(1\leq k,j\leq r).$

\begin{theorem}
\label{thm12}If $\mathcal{M}(n)\equiv\mathcal{M}^{N}(n)\geq0$ admits a flat
extension $\mathcal{M}(n+1)$, then $\mathcal{V}:=\mathcal{V}(\mathcal{M}%
(n+1))$ satisfies $\operatorname*{card}\mathcal{V}=r\;(\equiv
\operatorname*{rank}\mathcal{M}(n))$, and $\mathcal{V}\equiv\{v_{j}%
\}_{j=1}^{r}$ forms the support of the unique representing measure $\mu$ for
$\mathcal{M}(n+1)$. \ If $\mathcal{B}\equiv\{T^{i_{k}}\}_{k=1}^{r}$ is a
maximal linearly independent subset of columns of $\mathcal{M}(n)$, then
$W_{\mathcal{B},\mathcal{V}}$ is invertible, and $\mu=\sum_{i=1}^{r}\rho
_{j}\delta_{v_{j}}$, where $\rho\equiv(\rho_{1},...,\rho_{r})$ is uniquely
determined by $\rho^{t}=W_{\mathcal{B},\mathcal{V}}^{-1}(\beta_{i_{1}%
},...,\beta_{i_{r}})^{t}$. $\ $
\end{theorem}

Theorem \ref{thm12} describes $\mu$ in terms of $\mathcal{V}(\mathcal{M}%
(n+1))$. \ To compute the variety of any moment matrix $\mathcal{M}(n)$, we
may rely on the following general result. \ Given $n\geq1$, write $J\equiv
J(n):=\{j\in\mathbb{Z}_{+}^{N}:\left|  j\right|  \leq n\}$. \ Clearly, size
$\mathcal{M}(n)=\operatorname*{card}J(n)=\binom{N+n}{N}=\dim\{p\in
\mathbb{R}^{N}[t]:\deg p\leq n\}$.

\begin{proposition}
\label{prop1.3}Let $\mathcal{M}(n)\equiv\mathcal{M}^{N}(n)$ be a real moment
matrix, with columns $T^{j}$ indexed by $j\in J$, let $r:=\operatorname*{rank}%
\mathcal{M}(n)$, and let $\mathcal{B}\equiv\{T^{i}\}_{i\in I}$ be a maximal
linearly independent set of columns of $\mathcal{M}(n)$, where $I\subseteq J$
satisfies $\operatorname*{card}I=r$. \ For each index $j\in J\;\backslash\;I$,
let $q_{j}$ denote the unique polynomial in lin.span $\{t^{i}\}_{i\in I}$ such
that $T^{j}=q_{j}(T)$, and let $r_{j}(t):=t^{j}-q_{j}(t)$. \ Then
$\mathcal{V}(\mathcal{M}(n))$ is precisely the set of common zeros of
$\{r_{j}\}_{j\in J\;\backslash\;I}$.
\end{proposition}

Cases are known where $\beta^{(2n)}$ has no $\operatorname*{rank}%
\mathcal{M}(n)$-atomic $K_{\mathcal{Q}}$-representing measure, but does have a
\textit{finitely atomic} $K_{\mathcal{Q}}$-representing measure (cf.
\cite{tcmp6}, \cite{tcmp9}, \cite{Fia4}). \ It follows from Theorem
\ref{IntroTheorem1} that $\beta^{(2n)}$ has a finitely atomic representing
measure supported in $K_{\mathcal{Q}}$ if and only if $\mathcal{M}(n)(\beta)$
admits some positive moment matrix extension $\mathcal{M}(n+j)$, which in turn
admits a flat extension $\mathcal{M}(n+j+1)$ for which the unique successive
flat extensions $\mathcal{M}(n+j+k)$ satisfy $\mathcal{M}_{q_{i}}%
(n+j+k_{i})\geq0\;\ (1\leq i\leq m)$. \ We may estimate the minimum size of
$j$ as follows.\newpage

\begin{corollary}
\label{cor1.4}The $N$-dimensional real sequence $\beta^{(2n)}$ has a finitely
atomic representing measure supported in $K_{\mathcal{Q}}$ if and only if
$\mathcal{M}(n)(\beta)$ admits some positive moment matrix extension
$\mathcal{M}(n+j)$, with $j\leq2\binom{2n+N}{N}-n$, which in turn admits a
flat extension $\mathcal{M}(n+j+1)$ for which $\mathcal{M}_{q_{i}}%
(n+j+k_{i})\geq0\;\ (1\leq i\leq m)$.
\end{corollary}

If the conditions of Corollary \ref{cor1.4} hold, then the atoms and densities
of a finitely atomic $K_{\mathcal{Q}}$-representing measure for $\beta$ may be
computed by applying Theorem \ref{thm12} and Proposition \ref{prop1.3} to the
flat extension $\mathcal{M}(n+j+1)$. \ It is an open problem whether the
existence of a representing measure $\mu$ for $\beta^{(2n)}$ implies the
existence of a finitely atomic representing measure; such is the case, for
example, if $\mu$ has convergent moments of degree $2n+1$ (cf. \cite[Theorem
1.4]{tcmp8}, \cite{Put3}, \cite{Tch}).

We view Theorem \ref{IntroTheorem1} as our main result concerning existence of
minimal $K_{\mathcal{Q}}$-representing measures, and Theorem \ref{thm12}
primarily as a tool for computing such measures (cf. Example
\ref{IntroExample1} below). \ Note that Theorem \ref{thm12} applies to
arbitrary moment problems, not just the $K$-moment problem. \ Although Theorem
\ref{thm12} can also be regarded as an existence result, it may be very
difficult to utilize it in this way in specific examples. \ To explain this
viewpoint, we recall a result of \cite{EFP}. \ Let $\omega$ denote the
restriction of planar Lebesgue measure to the closed unit disk $\mathbb{\bar
{D}}$ and consider $\beta\equiv\beta^{(6)}[\omega]$ and $\mathcal{M\equiv
M}(3)(\beta)$; then $\operatorname*{rank}\mathcal{M}=10$. \ Flat extensions
$\mathcal{M}(4)$ of $\mathcal{M}$ exist in abundance and correspond to
$10$-atomic (minimal) cubature rules $\nu$ of degree $6$ for $\omega$. \ In
\cite{EFP} it is proved that no such rule $\nu$ is ``inside,'' i.e., with
$\operatorname*{supp}\nu\subseteq\mathbb{\bar{D}}$. \ The proof in \cite{EFP}
first characterizes the flat extensions $\mathcal{M}(4)$ in terms of algebraic
relations among the ``new moments'' of degree $7$ that appear in such
extensions. \ These relations lead to inequalities which ultimately imply
that, in Theorem \ref{IntroTheorem1}, $\mathcal{M}_{p}(4)$ cannot be positive
semi-definite, where $p(x,y):=1-x^{2}-y^{2}$. \ One could also try to
establish the nonexistence of $10$-atomic inside rules directly from Theorem
\ref{thm12}, without recourse to Theorem \ref{IntroTheorem1}. \ In this
approach one would first compute general formulas for the new moments of
degree $7$ in a flat extension $\mathcal{M}(4)$, use these moments to compute
the general form of $\mathcal{V}(\mathcal{M}(4))$, and then\ show that
$\mathcal{V}(\mathcal{M}(4))$ cannot be contained in $\mathbb{\bar{D}}$. \ As
a practical matter, however, this plan cannot be carried out; the new moments
comprise the solution of a system of $6$ quadratic equations in $8$ real
variables, and at present a program such as \textit{Mathematica} seems unable
to solve this system in a tractable form. \ For a problem such as this,
Theorem \ref{IntroTheorem1} seems indispensable. \ We illustrate the interplay
between Theorem \ref{IntroTheorem1}, Theorem \ref{thm12} and Proposition
\ref{prop1.3} in Example \ref{IntroExample1} below.

For measures in the plane ($N=2$), Theorem \ref{IntroTheorem1} is equivalent
to \cite[Theorem 1.6]{tcmp4}, which characterizes the existence of minimal
$K$-representing measures in the semi-algebraic case of the \textit{truncated
complex }$K$\textit{-moment problem} (with moments relative to monomials of
the form $\bar{z}^{i}z^{j}$). \ In \cite{tcmp4} we remarked that \cite[Theorem
1.6]{tcmp4} extended to truncated moment problems in any number of real or
complex variables. \ In \cite{Las1}, Lasserre developed applications of
\cite[Theorem 1.6]{tcmp4} to optimization problems in the plane. \ These
applications also extend to $\mathbb{R}^{N}$ $(N>2)$ (cf.\ \cite{Las1},
\cite{Las2}, \cite{Las3}), but they require the above mentioned generalization
of \cite[Theorem 1.6]{tcmp4} that we provide in Theorem \ref{IntroTheorem1}.
\ Lasserre's work motivated us to revisit our assertion in \cite{tcmp4}; we
then realized that there were unforeseen difficulties with the generalization,
particularly for the case when $N$ is odd. \ The purpose of Theorem
\ref{IntroTheorem1} is to provide the desired generalization.

The proofs of Theorem \ref{IntroTheorem1} and Corollary \ref{cor1.4} appear in
Section \ref{Exi}. \ In Theorem \ref{ThmExi.1} we characterize the existence
of minimal $K$-representing measures in the semi-algebraic case of the
truncated complex $K$-moment problem for measures on $\mathbb{C}^{m}$. \ The
equivalence of this result to the ``even'' case of Theorem \ref{IntroTheorem1}
($N=2d)$ is given in the first part of the proof of Theorem \ref{ThmExi.2};
this is based on the equivalence of the truncated moment problem for
$\mathbb{C}^{d}$ with the truncated real moment problem for $\mathbb{R}^{2d}$
(cf. Propositions \ref{prop12}, \ref{prop13}, \ref{PropositionEquiv} and
\ref{prop15}). \ The proof of Theorem \ref{IntroTheorem1} for $N=2d-1$, given
in the second part of the proof of Theorem \ref{ThmExi.2}, requires an
additional argument, based on the equivalence of a truncated moment problem
for $\mathbb{R}^{2d-1}$ with an associated moment problem for $\mathbb{R}%
^{2d}$.

We prove Theorem \ref{thm12} and Proposition \ref{prop1.3} in Section
\ref{Mom}. \ Theorem \ref{thm12} is new even for $N=2$. \ Previously, for
$N=2$ we knew that the measure $\mu$ of Theorems \ref{IntroTheorem1} and
\ref{thm12} could be computed with $\operatorname*{supp}\mu=\mathcal{V}%
(\mathcal{M}(r))$ \cite[p. 33]{tcmp1}, where $r:=\operatorname*{rank}%
\mathcal{M}(n)$ satisfies $r\leq\frac{(n+1)(n+2)}{2}$; but for $r>n+1$ this
entails iteratively generating the extensions $\mathcal{M}%
(n+2),...,\mathcal{M}(r)$. \ For $N>2$, we previously had no method for
computing $\mu$. \ In order to prove Theorem \ref{thm12} we first obtain some
results concerning truncated complex moment problems on $\mathbb{C}^{d}$.
\ Let $M(n)\equiv M^{d}(n)(\gamma)$ denote the moment matrix for a
$d$-dimensional complex multisequence $\gamma$ of degree $2n$, and let
$\mathcal{V}(M(n))$ denote the corresponding algebraic variety. \ In Theorem
\ref{propflat} we prove that if $M(n)\geq0$ admits a flat extension $M(n+1)$,
then the unique successive flat moment matrix extensions $\mathcal{M}%
(n+2),\mathcal{M}(n+3),...$ (cf. Theorem \ref{MomThm2}) satisfy $\mathcal{V}%
(M(n+1))=\mathcal{V}(M(n+2))=...$ . \ This result is used to prove Theorem
\ref{thmB}, which is the analogue of Theorem \ref{thm12} for the complex
moment problem. \ The proof of Theorem \ref{thm12} is then given in Theorem
\ref{thmA}, using Theorem \ref{thmB} and the ``equivalence'' results cited above.

In Section \ref{Loc} we study the \textit{localizing matrix} $M_{p}^{d}(n)$
corresponding to a complex moment matrix $M^{d}(n)$ and a polynomial
$p\in\mathbb{C}_{2n}^{d}[z,\bar{z}]$; \ Theorem \ref{ThmLoc.2} provides a
computational formula for $M_{p}^{d}(n)$ as a linear combination of certain
compressions of $M^{d}(n)$ corresponding to the monomial terms of $p$; an
analogous formula holds as well for \textit{real localizing matrices} (cf.
Theorem \ref{thm3.6}). \ In Section\ \ref{Fla}, we show that a flat extension
$M^{d}(n+1)$ of $M^{d}(n)\geq0$ induces flat extensions of positive localizing
matrices. \ Indeed, the flat extension $M^{d}(n+1)$ has unique successive flat
extensions $M^{d}(n+2),M^{d}(n+3),...$, and in Theorem \ref{ThmFla.1}, for
$p\in\mathbb{C}^{d}[z,\bar{z}]$, $\deg p=2k$ or $2k-1$, we prove that if
$M_{p}^{d}(n+k)\geq0$, then $M_{p}^{d}(n+k+1)$ is a flat, positive extension
of $M_{p}^{d}(n+k)$. \ In proving Theorem \ref{ThmFla.1} we follow the same
general plan as in the proof of \cite[Theorem 1.6]{tcmp4} (for moment problems
on $\mathbb{C}$), but we have streamlined the argument somewhat, placing more
emphasis on the abstract properties of flat extensions and less emphasis on
detailed calculations of the extensions; such calculations unnecessarily
complicated the argument given in \cite{tcmp4}. \ Theorem \ref{ThmFla.1} is
the main technical result that we need to prove Theorem \ref{IntroTheorem1}.

In the following example, we show the interaction of Theorem
\ref{IntroTheorem1}, Theorem \ref{thm12} and Proposition \ref{prop1.3} in a
$3$-dimensional cubature problem.

\begin{example}
\label{IntroExample1}We consider the cubature problem of degree $2$ for volume
measure $\mu\equiv\mu_{\mathbb{B}^{3}}$ on the closed unit ball $\mathbb{B}%
^{3}$ in $\mathbb{R}^{3}$ (cf. \cite{Str}). \ Thus $\beta\equiv\beta
^{(2)}=\left\{  \beta_{(i,j,k)}\right\}  _{i,j,k\geq0,\;i+j+k\leq2}$, where
$\beta_{(i,j,k)}:=\int_{\mathbb{B}^{3}}x^{i}y^{j}z^{k}\;d\mu$, i.e.,
$\beta_{(0,0,0)}=\frac{4\pi}{3}$, $\beta_{(1,0,0)}=\beta_{(0,1,0)}%
=\beta_{(0,0,1)}=0$, $\beta_{(2,0,0)}=\beta_{(0,2,0)}=\beta_{(0,0,2)}%
=\frac{4\pi}{15}$, $\beta_{(1,1,0)}=\beta_{(1,0,1)}=\beta_{(0,1,1)}=0$. \ The
moment matrix $\mathcal{M}^{3}(1)(\beta)$ has rows and columns indexed by $1$,
$X$, $Y$, $Z$; for $i\equiv(i_{1},i_{2},i_{3}),\;j\equiv(j_{1},j_{2},j_{3}%
)\in\mathbb{Z}_{+}^{3}$ with $\left|  i\right|  ,\left|  j\right|  \leq1$, the
entry in row $X^{i_{1}}Y^{i_{2}}Z^{i_{3}}$, column $X^{j_{1}}Y^{j_{2}}%
Z^{j_{3}}$, is $\beta_{(i_{1}+j_{1,\,}i_{2}+j_{2},\,i_{3}+j_{3})}$. Thus we
have $\mathcal{M}\equiv\mathcal{M}^{3}(1)(\beta)=$ diag $(\frac{4\pi}%
{3},\frac{4\pi}{15},\frac{4\pi}{15},\frac{4\pi}{15})$. $\ $We will use Theorem
\ref{IntroTheorem1} to construct a $\operatorname*{rank}\mathcal{M}$-atomic
representing measure for $\beta$ supported in $K=\mathbb{B}^{3}$. \ A moment
matrix extension $\mathcal{M}(2)$ of $\mathcal{M}$ admits a block
decomposition $\mathcal{M}(2)=\left(
\begin{array}
[c]{cc}%
\mathcal{M} & \mathcal{B}(2)\\
\mathcal{B}(2)^{t} & \mathcal{C}(2)
\end{array}
\right)  $, where $\mathcal{B}(2)$ includes ``new moments'' of degree $3$ and
$\mathcal{C}(2)$ is a moment matrix block of degree $4$; the rows and columns
of $M(2)$ are indexed by $1,$ $X,$ $Y,$ $Z,$ $X^{2},$ $YX,$ $ZX,$ $Y^{2},$
$ZY,$ $Z^{2}$ (see Section \ref{Mom} below). Clearly, $\mathcal{M}$ is
positive definite and invertible, so a flat extension $\mathcal{M}(2)$ is
determined by a choice of moments of degree $3$ such that $\mathcal{B}%
(2)^{t}\mathcal{M}^{-1}\mathcal{B}(2)$ has the form of a moment matrix block
$\mathcal{C}(2)$ (cf. the remarks following Theorem \ref{thmB}). \ Due to its
complexity, we are unable to compute the general solution $\mathcal{B}(2)$ to%
\begin{equation}
\mathcal{C}(2)=\mathcal{B}(2)^{t}\mathcal{M}^{-1}\mathcal{B}(2).
\label{IntroEq3}%
\end{equation}
Instead, we specify certain moments of degree $3$ as follows:%
\begin{align}
\beta_{(2,0,1)}  &  =\beta_{(2,1,0)}=\beta_{(1,1,1)}=\beta_{(0,2,1)}%
=\beta_{(0,1,2)}=0,\label{IntroEq4}\\
\beta_{(3,0,0)}  &  =\frac{1125\beta_{(1,2,0)}^{2}-16\pi^{2}}{1125\beta
_{(1,2,0)}},\;\;\;\beta_{(1,0,2)}=-\frac{16\pi^{2}}{1125\beta_{(1,2,0)}%
}\text{.}\nonumber
\end{align}
(Observe that we have left $\beta_{(1,2,0)}$, $\beta_{(0,3,0)}$ and
$\beta_{(0,0,3)}$ free.) \ With these choices, $\mathcal{B}(2)^{t}%
\mathcal{M}^{-1}\mathcal{B}(2)$ is a moment matrix block of degree $4$, and
$\mathcal{M}(2)\equiv\mathcal{M}(2)\{\beta_{(1,2,0)},\beta_{(0,3,0)}%
,\beta_{(0,0,3)}\}$ (defined by (\ref{IntroEq3})) is a flat extension of
$\mathcal{M}$. \ To show that $\beta$ admits a $4$-atomic $K$-representing
measure, we consider $p(x,y,z)=1-(x^{2}+y^{2}+z^{2})$, so that $K=K_{p}$
(where by $K_{p}$ we mean $K_{\mathcal{Q}}$ with $\mathcal{Q}\equiv\{p\}$).
\ Since $\deg p=2$, in Theorem \ref{IntroTheorem1} we have $n=k=1$; it thus
suffices to show that the flat extension $\mathcal{M}(2)$ corresponding to
(\ref{IntroEq4}) satisfies $\mathcal{M}_{p}(2)\geq0$. \ As we describe in
Section \ref{Loc} below, $\mathcal{M}_{p}(2)=\mathcal{M}_{1}(2)-(\mathcal{M}%
_{x^{2}}(2)+\mathcal{M}_{y^{2}}(2)+\mathcal{M}_{z^{2}}(2))$, where
$\mathcal{M}_{1}(2)=\mathcal{M},$ $\mathcal{M}_{x^{2}}(2)$ is the compression
of $\mathcal{M}(2)$ to rows and columns indexed by $X$, $X^{2}$, $YX$, $ZX,$
$\mathcal{M}_{y^{2}}(2)$ is the compression of $\mathcal{M}(2)$ to rows and
columns indexed by $Y$, $YX$, $Y^{2}$, $ZY$, and $\mathcal{M}_{z^{2}}(2)$ is
the compression of $\mathcal{M}(2)$ to rows and columns indexed by $Z$, $ZX$,
$ZY$, $Z^{2}$. \ From these observations, and using (\ref{IntroEq3}%
)--(\ref{IntroEq4}), it is straightforward to verify that
\[
\mathcal{M}_{p}(2)=\left(
\begin{array}
[c]{cccc}%
\frac{8\pi}{15} & -2\beta_{(3,0,0)} & -\beta_{(0,3,0)} & -\beta_{(0,0,3)}\\
&  &  & \\
-2\beta_{(3,0,0)} & f(\beta_{(1,2,0)}) & -\frac{15\beta_{(1,2,0)}%
\beta_{(0,3,0)}}{4\pi} & \frac{4\pi\beta_{(0,0,3)}}{75\beta_{(1,2,0)}}\\
&  &  & \\
-\beta_{(0,3,0)} & -\frac{15\beta_{(1,2,0)}\beta_{(0,3,0)}}{4\pi} &
g(\beta_{(1,2,0)},\beta_{(0,3,0)}) & 0\\
&  &  & \\
-\beta_{(0,0,3)} & \frac{4\pi\beta_{(0,0,3)}}{75\beta_{(1,2,0)}} & 0 &
h(\beta_{(1,2,0)},\beta_{(0,0,3)})
\end{array}
\right)  ,
\]
where
\begin{align*}
f(r)  &  :=-\frac{(1125r^{2}-300\pi r+16\pi^{2})(1125r^{2}+300\pi r+16\pi
^{2})}{168750\pi r^{2}},\\
g(r,s)  &  :=-\frac{2250r^{2}+1125s^{2}-64\pi^{2}}{300\pi},\\
h(r,t)  &  :=-\frac{-72000\pi^{2}r^{2}+512\pi^{4}+1265625r^{2}t^{2}}{337500\pi
r^{2}},
\end{align*}
and that $\mathcal{M}_{p}(2)$ is positive semi-definite if $\beta
_{(0,3,0)}=\beta_{(0,0,3)}=0$ and $\frac{2}{15}\sqrt{\frac{2}{5}}\pi\leq
\beta_{(1,2,0)}\leq\frac{4}{15}\sqrt{\frac{2}{5}}\pi$. \ Under these
conditions, Theorem \ref{IntroTheorem1} now implies the existence of a unique
$4$-atomic $($minimal$)$ representing measure $\zeta$ for $\mathcal{M}(2)$,
each of whose atoms lies in the closed unit ball. \ Theorem \ref{thm12}
implies that $\operatorname*{supp}\zeta=\mathcal{V}:=\mathcal{V}%
(\mathcal{M}(2))$. \ To compute the atoms of $\zeta$ via Proposition
\ref{prop1.3}, observe that in the column space of $\mathcal{M}(2)$ we have
the following linear dependence relations: $X^{2}=\frac{1}{5}1+\frac{1125\beta
_{(1,2,0)}^{2}-16\pi^{2}}{300\pi\beta_{(1,2,0)}}X,$ $XY=\frac{15\beta
_{(1,2,0)}}{4\pi}Y,$ $XZ=-\frac{4\pi}{75\beta_{(1,2,0)}}Z,$ $Y^{2}=\frac{1}%
{5}1+\frac{15\beta_{(1,2,0)}}{4\pi}X,$ $YZ=0,$ and $Z^{2}=\frac{1}%
{5}1-\frac{4\pi}{75\beta_{(1,2,0)}}X$; thus, $\mathcal{V}$ is determined by
the polynomials corresponding to these relations. \ A calculation shows that
$\mathcal{V}=\left\{  P_{i}\right\}  _{i=0}^{3}$, where $P_{i}\equiv
(x_{i},y_{i},z_{i})$ satisfies
\begin{align*}
P_{0}  &  =\left(  \frac{15\beta_{(1,2,0)}}{4\pi},-\frac{s}{4\sqrt{5}\pi
},0\right)  ,\;P_{1}=\left(  \frac{15\beta_{(1,2,0)}}{4\pi},\frac{s}{4\sqrt
{5}\pi},0\right)  ,\\
P_{2}  &  =\left(  -\frac{4\pi}{75\beta_{(1,2,0)}},0,-\frac{s}{75\beta
_{(1,2,0)}}\right)  ,\;P_{3}=\left(  -\frac{4\pi}{75\beta_{(1,2,0)}%
},0,\frac{s}{75\beta_{(1,2,0)}}\right)  ,
\end{align*}
with $s:=\sqrt{1125\beta_{(1,2,0)}^{2}+16\pi^{2}}$. \ The measure $\zeta$ is
thus of the form $\zeta=\sum_{i=0}^{3}\rho_{_{i}}\delta_{P_{_{i}}}$. \ To
compute the densities $\rho_{i}$ using Theorem \ref{thm12}, consider the basis
$\mathcal{B}:=\{1,X,Y,Z\}$ for $\mathcal{C}_{\mathcal{M}(1)}$ and let
\[
W=\left(
\begin{array}
[c]{cccc}%
1 & 1 & 1 & 1\\
x_{0} & x_{1} & x_{2} & x_{3}\\
y_{0} & y_{1} & y_{2} & y_{3}\\
z_{0} & z_{1} & z_{2} & z_{3}%
\end{array}
\right)  .
\]
Following Theorem \ref{thm12}, $\rho\equiv\left(  \rho_{0},\rho_{1},\rho
_{2},\rho_{3}\right)  $ is uniquely determined by
\[
\rho^{t}=W^{-1}\left(  \beta_{(0,0,0)},\beta_{(1,0,0)},\beta_{(0,1,0)}%
,\beta_{(0,0,1)}\right)  ^{t}=W^{-1}\left(  \frac{4\pi}{3},0,0,0\right)
^{t},
\]
and thus%
\[
\rho_{0}=\rho_{1}=\frac{32\pi^{3}}{3(1125\beta_{(1,2,0)}^{2}+16\pi^{2}%
)},\;\rho_{2}=\rho_{3}=\frac{750\beta_{(1,2,0)}^{2}\pi}{1125\beta
_{(1,2,0)}^{2}+16\pi^{2}}.
\]
For a concrete numerical example, we can take $\beta_{(1,2,0)}=\frac{4}%
{15}\sqrt{\frac{2}{5}}\pi$, and obtain $\rho_{0}=\rho_{1}=\frac{2}{9}%
\pi,\;\rho_{2}=\rho_{3}=\frac{4}{9}\pi$, and $P_{0}=\left(  \sqrt{\frac{2}{5}%
},-\sqrt{\frac{3}{5}},0\right)  ,\;P_{1}=\left(  \sqrt{\frac{2}{5}}%
,\sqrt{\frac{3}{5}},0\right)  ,$ $P_{2}=\left(  -\sqrt{\frac{1}{10}}%
,0,-\sqrt{\frac{3}{10}}\right)  ,$ and $P_{3}=\left(  -\sqrt{\frac{1}{10}%
},0,\sqrt{\frac{3}{10}}\right)  $. $\ $Note that
\[
\mathcal{M}_{p}(2)=\left(
\begin{array}
[c]{cc}%
\frac{8\pi}{15} & -\frac{4}{15}\sqrt{\frac{2}{5}}\pi\\
-\frac{4}{15}\sqrt{\frac{2}{5}}\pi & \frac{4\pi}{75}%
\end{array}
\right)
{\textstyle\bigoplus}
(0)%
{\textstyle\bigoplus}
(\frac{4\pi}{25}),
\]
so $\operatorname*{rank}\mathcal{M}(1)-\operatorname*{rank}\mathcal{M}%
_{p}(2)=2$, and (as Theorem \ref{IntroTheorem1} predicts) there are two
points, $P_{0}$ and $P_{1}$, that lie on the unit sphere. \
\endproof
\end{example}

We pause to locate Theorem \ref{IntroTheorem1} within the extensive literature
on the $K$-moment problem (cf. \cite{Akh}, \cite{BeCJ}, \cite{BeMa},
\cite{Fug}, \cite{KrNu}, \cite{Rez}, \cite{ShTa}, \cite{StSz2}). \ A classical
theorem of M. Riesz \cite[Section 5]{Rie} provides a solution to the full
$K$-moment problem on $\mathbb{R}$, as follows. Given a real sequence
$\beta\equiv\left\{  \beta_{i}\right\}  _{i=0}^{\infty}$ and a closed set
$K\subseteq\mathbb{R}$, there exists a positive Borel measure $\mu$ on
$\mathbb{R}$ such that $\beta_{i}=\int t^{i}\,d\mu$ \ ($i\geq0$) and
$\operatorname*{supp}\mu\subseteq K$ if and only if each polynomial
$p\in\mathbb{C}\left[  t\right]  ,\;p(t)=\sum_{i=0}^{N}a_{i}t^{i}$, with
$p|_{K}\geq0$, satisfies $\sum_{i=0}^{N}a_{i}\beta_{i}\geq0$. For a general
closed set $K\subseteq\mathbb{R}$ there is no concrete description of the case
$p|_{K}\geq0$, so it may be very difficult to verify the Riesz hypothesis for
a particular $\beta$.

In \cite{Hav}, Haviland extended Riesz's theorem to $\mathbb{R}^{N}$ $(N>1)$
and also showed that for several semi-algebraic sets $K$, the Riesz hypothesis
can be checked by concrete positivity tests. \ Indeed, by combining the
generalized Riesz hypothesis with concrete descriptions of non-negative
polynomials on $\mathbb{R}$, $\left[  0,+\infty\right]  ,\left[  a,b\right]
$, or the unit circle, Haviland recovered classical solutions to the full
moment problems of Hamburger, Stieltjes, Hausdorff, and Herglotz \cite{Hav}.
\ More recently, for the case of the closed unit disk, Atzmon \cite{Atz} found
a concrete solution to the full $K$-moment problem using subnormal operator
theory, and Putinar \cite{Put1} subsequently presented an alternate solution
using hyponormal operator theory.

In \cite{Cas}, Cassier initiated the study of the $K$-moment problem for
compact subsets of $\mathbb{R}^{N}$. For the case when $K$ is compact and
semi-algebraic, Schm\"{u}dgen \cite{Sch} used real algebraic geometry to solve
the full $K$-moment problem in terms of concrete positivity tests. Using
infinite moment matrices, we may paraphrase Schm\"{u}dgen's theorem as
follows: a full multi-sequence $\beta\equiv\beta^{(\infty)}=\left\{  \beta
_{i}\right\}  _{^{i\in\mathbb{Z}_{+}^{N}}}$ has a representing measure
supported in a compact semi-algebraic set $K_{\mathcal{Q}}$ if and only if
$\mathcal{M}^{N}(\infty)(\beta)\geq0$ and $\mathcal{M}_{q}^{N}(\infty
)(\beta)\geq0$ for every polynomial $q$ that is a product of distinct $q_{i}$.
Schm\"{u}dgen's approach, using real algebra, is to concretely describe the
polynomials nonnegative on $K_{\mathcal{Q}}$ (as above) and to then apply the
Riesz-Haviland criterion. \ Putinar and Vasilescu \cite{PuVa1} subsequently
provided a reduced set of testing polynomials $q$ (see also \cite{Dem}).
\ Recently, Powers and Scheiderer \cite{PoSc} characterized the non-compact
semi-algebraic sets $K_{\mathcal{Q}}$ for which a generalized
Schm\"{u}dgen-type theorem is valid. \ Indeed, recent advances in real algebra
make it possible to concretely describe the polynomials nonnegative on certain
noncompact semi-algebraic sets (\cite{KuMa}, \cite{PoRe1}, \cite{PoRe2},
\cite{PoSc}, \cite{Pre}, \cite{Put2}, \cite{Sche}), so as to establish moment
theorems via the previously intractable Riesz-Haviland approach.

There is at present no viable analogue of the Riesz-Haviland criterion for
truncated moment problems. \ Theorem \ref{IntroTheorem1} is motivated by the
above results for the full $K$-moment problem and also by a recent result of
J. Stochel \cite{Sto2} which shows that the truncated $K$-moment problem is
actually more general than the full $K$-moment problem. Stochel's result in
\cite{Sto2} is stated for the complex multidimensional moment problem, but we
may paraphrase it for the real moment problem as follows.

\begin{theorem}
\label{IntroTheorem2} (cf.\ \cite{Sto2}). Let $K$ be a closed subset of
$\mathbb{R}^{N}$ $(N>1)$. A real multisequence $\beta\equiv\beta^{(\infty
)}=\left\{  \beta_{i}\right\}  _{i\in\mathbb{Z}_{+}^{N}}$ has a $K$%
-representing measure if, and only if, for each $n>0$, $\beta^{(2n)}%
\equiv\left\{  \beta_{i}\right\}  _{i\in\mathbb{Z}_{+}^{N},\;\left|  i\right|
\leq2n}$ has a $K$-representing measure.
\end{theorem}

For the semi-algebraic case $(K=K_{\mathcal{Q}})$, Theorem \ref{IntroTheorem1}
addresses the existence of finitely atomic $K$-representing measures for
$\beta^{(2n)}$ with the fewest atoms possible. Concerning the existence of a
flat extension $\mathcal{M}^{N}(n+1)$ in Theorem \ref{IntroTheorem1}, there is
at present no satisfactory general test available, so in this sense Theorem
\ref{IntroTheorem1} is ``abstract.'' However, in certain special cases,
concrete solutions to the flat extension problem have been found
(\cite{tcmp2}, \cite{tcmp3}). \ For example, consider the case of the
parabolic moment problem, where $q(x,y)=0$ represents a parabola in
$\mathbb{R}^{2}$. \ Theorem \ref{IntroTheorem1} implies that $\beta^{(2n)}$
has a $\operatorname*{rank}\mathcal{M}^{2}(n)$-atomic representing measure
supported in $\mathcal{Z}(q)$ if and only if $\mathcal{M}^{2}(n)(\beta)$ is
positive and admits a flat extension $\mathcal{M}^{2}(n+1)$ for which
$\mathcal{M}_{q}^{2}(n+1)=0$. \ In \cite{tcmp7} we obtained the following
concrete characterization of this case.

\begin{theorem}
\label{IntroTheorem4}(\cite[Theorem 2.2]{tcmp7}) Let $q(x,y)=0$ denote a
parabola in $\mathbb{R}^{2}$. \ The following statements are equivalent for
$\beta\equiv\beta^{(2n)}$:\newline (i) $\beta$ has a representing measure
supported in $\mathcal{Z}(q)$;\newline (ii) $\beta$ has a (minimal)
$\operatorname*{rank}\mathcal{M}^{2}(n)(\beta)$-atomic representing measure
supported in $\mathcal{Z}(q)$ (cf. Theorem\ \ref{IntroTheorem1});\newline
(iii) $\mathcal{M}^{2}(n)(\beta)$ is positive and \textit{recursively
generated} (cf. Section \ref{Mom}), there is a column dependence relation
$q(X,Y)=0$, and $\operatorname*{card}\mathcal{V}(\mathcal{M}^{2}%
(n)(\beta))\geq\operatorname*{rank}\mathcal{M}^{2}(n)(\beta)$. \ 
\end{theorem}%

\noindent
Analogues of Theorem \ref{IntroTheorem4} for all other curves of degree $2$
appear in \cite{tcmp5}, \cite{tcmp6}, \cite{tcmp9}, \cite{Fia4}. \ The full
moment problem on a curve of degree $2$ had previously been concretely solved
in \cite{Sto1} (cf. \cite{StSz1}); an alternate solution appears in
\cite{PoSc}.

\textit{Acknowledgments}. Example \ref{IntroExample1} was obtained using
calculations with the software tool \textit{Mathematica \cite{Wol}}. \ The
authors are grateful to the referee for several suggestions that helped
improve the presentation.

\section{\label{Mom}Moment matrices}

Let $\mathbb{C}_{r}^{d}\left[  z,\bar{z}\right]  $ denote the space of
polynomials with complex coefficients in the indeterminates $z\equiv
(z_{1},...,z_{d})$ and $\bar{z}\equiv(\bar{z}_{1},...,\bar{z}_{d})$, with
total degree at most $r$; \ thus $\dim\mathbb{C}_{r}^{d}\left[  z,\bar
{z}\right]  =\eta(d,r):=\left(
\begin{array}
[c]{c}%
r+2d\\
2d
\end{array}
\right)  $. \ For $i\equiv(i_{1},...,i_{d})\in\mathbb{Z}_{+}^{d}$, let
$\left|  i\right|  :=i_{1}+...+i_{d}$ and let $z^{i}:=z_{1}^{i_{1}}\cdots
z_{d}^{i_{d}}$. \ Given a complex sequence $\gamma\equiv\gamma^{(s)}%
=\{\gamma_{ij}\}_{i,j\in\mathbb{Z}_{+}^{d}}$, $\left|  i\right|  +\left|
j\right|  \leq s$, the \textit{truncated complex moment problem} for $\gamma$
entails determining necessary and sufficient conditions for the existence of a
positive Borel measure $\nu$ on $\mathbb{C}^{d}$ such that
\begin{equation}
\gamma_{ij}=\int\bar{z}^{i}z^{j}\;d\nu\;\;(\equiv\int\bar{z}_{1}^{i_{1}}%
\cdots\bar{z}_{d}^{i_{d}}z_{1}^{j_{1}}\cdots z_{d}^{j_{d}}\;d\nu
(z_{1},...,z_{d},\bar{z}_{1},...,\bar{z}_{d}))\;\;(\left|  i\right|  +\left|
j\right|  \leq s).\ \label{eq.Mom1}%
\end{equation}
A measure $\nu$ as in (\ref{eq.Mom1}) is a \textit{representing measure} for
$\gamma^{(s)}$; \ if $K\subseteq\mathbb{C}^{d}$ is a closed set and
$\operatorname*{supp}\nu\subseteq K$, then $\nu$ is a $K$\textit{-representing
measure} for $\gamma^{(s)}$. \ 

In the sequel we focus on the case when $s$ is even, say $s=2n$. \ In this
case, the moment data $\gamma^{(2n)}$ determine the \textit{moment matrix}
$M(n)\equiv M^{d}(n)(\gamma)$ that we next describe. \ The size of $M(n)$ is
$\eta(d,n)$, with rows and columns $\{\bar{Z}^{i}Z^{j}\}_{i,j\in\mathbb{Z}%
_{+}^{d},\left|  i\right|  +\left|  j\right|  \leq n}$, indexed by the
lexicographic ordering of the monomials in $\mathbb{C}_{n}^{d}\left[
z,\bar{z}\right]  $; for $d=2,n=2$, this ordering is $1,Z_{1},Z_{2},\bar
{Z}_{1},\bar{Z}_{2},$ $Z_{1}^{2},Z_{1}Z_{2},\bar{Z}_{1}Z_{1},\bar{Z}_{2}%
Z_{1},Z_{2}^{2},$ $\bar{Z}_{1}Z_{2},\bar{Z}_{2}Z_{2},\bar{Z}_{1}^{2},\bar
{Z}_{1}\bar{Z}_{2},\bar{Z}_{2}^{2}$. \ The entry of $M(n)$ in row $\bar{Z}%
^{i}Z^{j}$, column $\bar{Z}^{k}Z^{\ell}$ is $\gamma_{i+\ell,k+j}\;(\left|
i\right|  +\left|  j\right|  ,\left|  k\right|  +\left|  \ell\right|  \leq
n)$. \ By a \textit{representing measure for }$M(n)$ we mean a representing
measure for $\gamma$.

For $p\in\mathbb{C}_{n}^{d}\left[  z,\bar{z}\right]  $, $p(z,\bar{z}%
)\equiv\sum_{r,s\in\mathbb{Z}_{+}^{d},\left|  r\right|  +\left|  s\right|
\leq n}a_{rs}\bar{z}^{r}z^{s}$, we set $\hat{p}:=(a_{rs})$; $\hat{p}$ is the
coefficient vector of $p$ relative to the basis for $\mathbb{C}_{n}^{d}\left[
z,\bar{z}\right]  $ consisting of the monomials $\{\bar{z}^{i}z^{j}%
\}_{i,j\in\mathbb{Z}_{+}^{d},\left|  i\right|  +\left|  j\right|  \leq n}$ in
lexicographic order. \ We recall the \textit{Riesz functional} $\Lambda
\equiv\Lambda_{\gamma}:\mathbb{C}_{2n}^{d}\left[  z,\bar{z}\right]
\rightarrow\mathbb{C}$, defined by $\Lambda(\sum_{r,s\in\mathbb{Z}_{+}%
^{d},\left|  r\right|  +\left|  s\right|  \leq2n}b_{rs}\bar{z}^{r}z^{s}%
):=\sum_{r,s\in\mathbb{Z}_{+}^{d},\left|  r\right|  +\left|  s\right|  \leq
2n}b_{rs}\gamma_{rs}$. \ The matrix $M^{d}(n)(\gamma)$ is uniquely determined
by
\begin{equation}
\left\langle M^{d}(n)(\gamma)\hat{f},\hat{g}\right\rangle :=\Lambda_{\gamma
}(f\bar{g}),f,g\in\mathbb{C}_{n}^{d}\left[  z,\bar{z}\right]
.\ \label{eq.Mom2}%
\end{equation}
If $\gamma$ has a representing measure $\nu$, then $\Lambda_{\gamma}(f\bar
{g})=\int f\bar{g}\;d\nu$; \ in particular, $\left\langle M^{d}(n)(\gamma
)\hat{f},\hat{f}\right\rangle =\int\left|  f\right|  ^{2}\;d\nu\geq0$, so
$M^{d}(n)(\gamma)$ is positive semidefinite in this case. \ 

Corresponding to $p\in\mathbb{C}_{n}^{d}\left[  z,\bar{z}\right]  $,
$p(z,\bar{z})\equiv\sum a_{rs}\bar{z}^{r}z^{s}$ (as above), we may define an
element in $\mathcal{C}_{M(n)}$, the column space of $M(n)$, by $p(Z,\bar
{Z}):=\sum a_{rs}\bar{Z}^{r}Z^{s}$; the following result will be used in the
sequel to locate the support of a representing measure.

\begin{proposition}
\label{MomProp1}(\cite[(7.4)]{tcmp1}) Suppose $\nu$ is a representing measure
for $\gamma^{(2n)}$, let $p\in\mathbb{C}_{n}^{d}\left[  z,\bar{z}\right]  $,
and let $\mathcal{Z}(p):=\{z\in\mathbb{C}^{d}:p(z,\bar{z})=0\}$. \ Then
$\operatorname*{supp}\nu\subseteq\mathcal{Z}(p)$ if and only if $p(Z,\bar
{Z})=0$.
\end{proposition}

It follows from Proposition \ref{MomProp1} that if $\gamma^{(2n)}$ has a
representing measure, then $M^{d}(n)(\gamma)$ is \textit{recursively
generated} in the following sense:
\begin{equation}
p,q,pq\in\mathbb{C}_{n}^{d}\left[  z,\bar{z}\right]  ,p(Z,\bar{Z}%
)=0\Rightarrow(pq)(Z,\bar{Z})=0. \label{eq.Mom3}%
\end{equation}

We define the \textit{variety} of $M(n)$ by $\mathcal{V}(M(n)):=\bigcap
_{p\in\mathbb{C}_{n}^{d}\left[  z,\bar{z}\right]  ,p(Z,\bar{Z})=0}%
\mathcal{Z}(p)$; we sometimes refer to $\mathcal{V}(M(n))$ as $\mathcal{V}%
(\gamma)$. \ Proposition \ref{MomProp1} implies that if $\nu$ is a
representing measure for $\gamma^{(2n)}$, then $\operatorname*{supp}%
\nu\subseteq\mathcal{V}(\gamma)$ and, moreover, that
\begin{equation}
\operatorname*{card}\mathcal{V}(\gamma)\geq\operatorname*{card}%
\operatorname*{supp}\nu\geq\operatorname*{rank}M^{d}(n)(\gamma)\;\;\text{(cf.
\cite[(7.6)]{tcmp1}).} \label{eq.Mom4}%
\end{equation}

The following result characterizes the existence of ``minimal,'' i.e.,
$\operatorname*{rank}M(n)$-atomic, representing measures.

\begin{theorem}
\label{MomThm2}(\cite[Corollary 7.9 and Theorem 7.10]{tcmp1}) $\gamma^{(2n)}$
has a $\operatorname*{rank}M^{d}(n)(\gamma)$-atomic representing measure if
and only if $M(n)\equiv M^{d}(n)(\gamma)$ is positive semidefinite and $M(n)$
admits an extension to a moment matrix $M(n+1)\equiv M^{d}(n+1)(\tilde{\gamma
})$ satisfying $\operatorname*{rank}M(n+1)=\operatorname*{rank}M(n)$. \ In
this case, $M(n+1)$ admits unique successive rank-preserving positive moment
matrix extensions $M(n+2),M(n+3),...$, and there exists a
$\operatorname*{rank}M(n)$-atomic representing measure for $M(\infty)$. \ 
\end{theorem}

Various concrete sufficient conditions are known for the existence of the
rank-preserving extension $M(n+1)$ described in Theorem \ref{MomThm2},
particularly when $d=1$ (moment problems in the plane) \cite{tcmp1},
\cite{tcmp2}, \cite{tcmp3}, \cite{tcmp5}, \cite{tcmp6}, \cite{tcmp7},
\cite{tcmp9}; \ for general $d$, an important sufficient condition is that
$M^{d}(n)(\gamma)$ is positive semidefinite and \textit{flat}, i.e.,
$\operatorname*{rank}M^{d}(n)(\gamma)$\newline $=\operatorname*{rank}%
M^{d}(n-1)(\gamma)$\ \cite[Theorem 7.8]{tcmp1}.

We now present the complex version of Theorem \ref{thm12}.

\begin{theorem}
\label{thmB}If $M(n)\equiv M^{d}(n)\geq0$ admits a rank-preserving extension
$M(n+1)$, then $\mathcal{V}:=\mathcal{V}(M(n+1))$ satisfies
$\operatorname*{card}\mathcal{V}=r\;(\equiv\operatorname*{rank}M(n))$, and
$\mathcal{V}\equiv\{\omega_{j}\}_{j=1}^{r}$ forms the support of the unique
representing measure $\nu$ for $M(n+1)$. \ If $\mathcal{B}\equiv\{\bar
{Z}^{i_{k}}Z^{j_{k}}\}_{k=1}^{r}$ is a maximal linearly independent subset of
columns of $M(n)$, then the $r\times r$ matrix $W_{\mathcal{B},\mathcal{V}}$
(whose entry in row $m$, column $k$ is $\bar{\omega}_{k}^{i_{m}}\omega
_{k}^{j_{m}}$) is invertible, and $\nu=\sum_{j=1}^{r}\rho_{j}\delta
_{\omega_{j}}$, where $\rho\equiv(\rho_{1},...,\rho_{r})$ is uniquely
determined by $\rho^{t}=W_{\mathcal{B},\mathcal{V}}^{-1}(\gamma_{i_{1},j_{1}%
},...,\gamma_{i_{r},j_{r}})^{t}$.
\end{theorem}

Toward the proof of Theorem \ref{thmB}, we begin with some remarks concerning
positive matrix extensions. \ Let $\tilde{A}\equiv\left(
\begin{array}
[c]{cc}%
A & B\\
B^{\ast} & C
\end{array}
\right)  $ be a block matrix. \ A result of Smul'jan \cite{Smu} shows that
$\tilde{A}\geq0$ if and only if $A\geq0$ and there exists a matrix $W$ such
that $B=AW$ and $C\geq W^{\ast}AW$. \ In this case, $W^{\ast}AW$ is
independent of $W$ satisfying $B=AW$, and the matrix $[A;B]:=\left(
\begin{array}
[c]{cc}%
A & B\\
B^{\ast} & W^{\ast}AW
\end{array}
\right)  $ is positive and satisfies $\operatorname*{rank}%
[A;B]=\operatorname*{rank}A$; conversely, any rank-preserving positive
extension $\tilde{A}$ of $A$ is of this form. \ We refer to such a
rank-preserving extension as a \textit{flat extension} of $A$. \ Now, a moment
matrix $M^{d}(n+1)$ admits a block decomposition $M(n+1)=\left(
\begin{array}
[c]{cc}%
M(n) & B(n+1)\\
B(n+1)^{\ast} & C(n+1)
\end{array}
\right)  $; \ thus a positive moment matrix $M(n)$ admits a flat (positive)
moment matrix extension $M(n+1)$ if and only if there is a choice of moments
of degree $2n+1$ and a matrix $W$ such that $B(n+1)=M(n)W$ and $W^{\ast}M(n)W$
has the form of a moment matrix block $C(n+1)$, i.e., $\left[
M(n);B(n+1)\right]  $ is a moment matrix.

Consider again a positive extension $\tilde{A}$ of $A$ (as above). \ The
Extension Principle (\cite[Proposition 3.9]{tcmp1}, \cite[Proposition
2.4]{FiaCM}) implies that each linear dependence relation in the columns of
$A$ extends to the columns of $\left(
\begin{array}
[c]{c}%
A\\
B^{\ast}%
\end{array}
\right)  $ in $\tilde{A}$. \ In the case when $M(n+1)$ is a positive extension
of $M(n)$, it follows that $\mathcal{V}(M(n+1))\subseteq\mathcal{V}(M(n))$; we
will use this relation frequently in the sequel, without further reference.

Now recall from Theorem \ref{MomThm2} that if $M(n)\geq0$ admits a flat
extension $M(n+1)$, then $M(n+1)$ admits a unique flat extension $M(n+2)$.
\ Indeed, every column of $M(n+1)$ of total degree $n+1$ is a linear
combination of columns corresponding to monomials of total degree at most $n$;
we can write this as
\begin{equation}
\bar{Z}^{i}Z^{j}=p_{i,j}(Z,\bar{Z})\;\;(p_{i,j}\in\mathbb{C}_{n}^{d}\left[
z,\bar{z}\right]  ;\;\left|  i\right|  +\left|  j\right|  =n+1).
\label{flatness}%
\end{equation}
Then the unique flat extension $M(n+2)$ is given by
\begin{equation}
\bar{Z}^{i}Z^{j}=\left\{
\begin{array}
[c]{cc}%
(z_{\ell}p_{i,j-\varepsilon(\ell)})(Z,\bar{Z})\;\; & \text{if }j_{\ell}%
\geq1\text{ for some }\ell=1,...,d\\
(\bar{z}_{k}p_{i-\varepsilon(k),0})(Z,\bar{Z}) & \text{if }j=0\text{ and
}i_{k}\geq1\text{ for some }k=1,...,d
\end{array}
\right.  \label{zij}%
\end{equation}
$(\left|  i\right|  +\left|  j\right|  =n+2)$, where $\varepsilon
(\ell):=(0,...,0,\overset{\ell}{1},0,...,0)$. \ ($\bar{Z}^{i}Z^{j}$ is
independent of the choice of $j_{\ell}$ or $i_{k}$; cf. \cite[Theorem
7.8]{tcmp1}.)

Suppose $M(n)\geq0$ admits a flat extension $M(n+1)$; the following result
implies that the unique rank-preserving extensions $M(n+2),M(n+3),...,$ are
also variety-preserving; this is a key ingredient in the proof of Theorem
\ref{thm12} and may be of independent interest.

\begin{theorem}
\label{propflat}Assume that $M(n)\equiv M^{d}(n)\geq0$ admits a flat extension
$M(n+1)$. \ Then \newline $\mathcal{V}(M(n+2))=\mathcal{V}(M(n+1)).$
\end{theorem}

\begin{proof}
Recall that $\mathcal{V}(M(n+2))\subseteq\mathcal{V}(M(n+1))$; to prove the
reverse inclusion, it suffices to show that if $\omega\in\mathcal{V}(M(n+1))$,
and $f\in\mathbb{C}_{n+2}^{d}\left[  z,\bar{z}\right]  $ satisfies
$f(Z,\bar{Z})=0$ in $\mathcal{C}_{M(n+2)}$, then $f(\omega,\bar{\omega})=0$.
\ As discussed above, the flat extension $M(n+2)$ admits a decomposition%
\[
M(n+2)=\left(
\begin{array}
[c]{cc}%
M(n+1) & M(n+1)W\\
W^{\ast}M(n+1) & W^{\ast}M(n+1)W
\end{array}
\right)  .
\]
Write $f=g+h$, where $g\in\mathbb{C}_{n+1}^{d}\left[  z,\bar{z}\right]  $, and
$h(z,\bar{z})\equiv\sum_{\left|  i\right|  +\left|  j\right|  =n+2}h_{i,j}%
\bar{z}^{i}z^{j}$. \ Recall that $\hat{f}\in\mathbb{C}^{\eta(d,n+2)}$ and
$\hat{g}\in\mathbb{C}^{\eta(d,n+1)}$ denote the coefficient vectors of $f$ and
$g$ relative to the bases of monomials in lexicographic order. \ Let
$\tilde{h}\in\mathbb{C}^{\eta(d,n+2)-\eta(d,n+1)}$ denote the coefficient
vector of $h$ relative to the monomials of degree $n+2$ in lexicographic
order; thus $\hat{f}=\left(
\begin{array}
[c]{c}%
\hat{g}\\
\tilde{h}%
\end{array}
\right)  $. \ Now,
\[
f(Z,\bar{Z})=M(n+2)\hat{f}=\left(
\begin{array}
[c]{c}%
M(n+1)\hat{g}+M(n+1)W\tilde{h}\\
W^{\ast}M(n+1)\hat{g}+W^{\ast}M(n+1)W\tilde{h}%
\end{array}
\right)  ,
\]
so $f(Z,\bar{Z})=0$ implies
\begin{equation}
M(n+1)(\hat{g}+W\tilde{h})=0.\ \label{zero}%
\end{equation}
We seek to associate $\hat{g}+W\tilde{h}$ with the coefficient vector $\hat
{q}$ of some polynomial $q\in\mathbb{C}_{n+1}^{d}\left[  z,\bar{z}\right]  $,
and to this end we first describe an explicit formula for $W$.

Recall that $M(n+1)W=B(n+2)$, and that the columns of $B(n+2)$ are associated
with the monomials $\bar{z}^{i}z^{j}\;\;((\left|  i\right|  +\left|  j\right|
=n+2)$. \ For $(i,j)\in\mathbb{Z}_{+}^{d}\times\mathbb{Z}_{+}^{d}$ with
$\left|  i\right|  +\left|  j\right|  =n+2$, the $(i,j)$-th column of $B(n+2)$
is, on one hand $M(n+1)W\widetilde{\bar{z}^{i}z^{j}}$, while on the other hand
it equals $\left[  (z_{\ell}p_{i,j-\varepsilon(\ell)})(Z,\bar{Z})\right]
_{\eta(d,n+1)}$ or $\left[  (\bar{z}_{k}p_{i-\varepsilon(k),0})(Z,\bar
{Z})\right]  _{\eta(d,n+1)}$, by (\ref{zij}). \ Since the polynomials
$z_{\ell}p_{i,j-\varepsilon(\ell)}$ and $\bar{z}_{k}p_{i-\varepsilon(k),0}$
belong to $\mathbb{C}_{n+1}^{d}\left[  z,\bar{z}\right]  $, we can write
\[
\left[  (z_{\ell}p_{i,j-\varepsilon(\ell)})(Z,\bar{Z})\right]  _{\eta
(d,n+1)}=M(n+1)(z_{\ell}p_{i,j-\varepsilon(\ell)})\symbol{94}%
\]
and
\[
\left[  (\bar{z}_{k}p_{i-\varepsilon(k),0})(Z,\bar{Z})\right]  _{\eta
(d,n+1)}=M(n+1)(\bar{z}_{k}p_{i-\varepsilon(k),0})\symbol{94}.
\]
It follows at once that $W$ can be given by
\begin{align}
&  W\widetilde{\bar{z}^{i}z^{j}}\nonumber\\
&  =\left\{
\begin{array}
[c]{cc}%
(z_{\ell}p_{i,j-\varepsilon(\ell)})\symbol{94} & \text{if }j_{\ell}\geq1\text{
for some }\ell=1,...,d\\
(\bar{z}_{k}p_{i-\varepsilon(k),0})\symbol{94} & \text{if }j=0\text{ and
}i_{k}\geq1\text{ for some }k=1,...,d
\end{array}
\right.  \;(\left|  i\right|  +\left|  j\right|  =n+2). \label{eq28}%
\end{align}
We now consider $W\tilde{h}$. \ Since $\tilde{h}\equiv\sum_{\left|  i\right|
+\left|  j\right|  =n+2}h_{i,j}\widetilde{\bar{z}^{i}z^{j}}$, it follows from
(\ref{eq28}) that
\begin{align*}
W\tilde{h}  &  =\sum_{\left|  i\right|  +\left|  j\right|  =n+2,j\neq0}%
h_{i,j}(z_{\ell}p_{i,j-\varepsilon(\ell)})\symbol{94}+\sum_{\left|  i\right|
=n+2}h_{i,0}(\bar{z}_{k}p_{i-\varepsilon(k),0})\symbol{94}\\
&  =[\sum_{\left|  i\right|  +\left|  j\right|  =n+2,j\neq0}h_{i,j}z_{\ell
}p_{i,j-\varepsilon(\ell)}+\sum_{\left|  i\right|  =n+2}h_{i,0}\bar{z}%
_{k}p_{i-\varepsilon(k),0}]\symbol{94}.
\end{align*}
\ Now we set
\begin{gather*}
q(z,\bar{z}):=g(z,\bar{z})+\sum_{\left|  i\right|  +\left|  j\right|
=n+2,j\neq0}h_{i,j}(z_{\ell}p_{i,j-\varepsilon(\ell)})(z,\bar{z})\\
+\sum_{\left|  i\right|  =n+2}h_{i,0}(\bar{z}_{k}p_{i-\varepsilon
(k),0})(z,\bar{z})\in\mathbb{C}_{n+1}^{d}\left[  z,\bar{z}\right]  .
\end{gather*}

Observe that in $\mathcal{C}_{M(n+1)}$,
\begin{align*}
q(Z,\bar{Z})  &  =M(n+1)\hat{q}=M(n+1)\hat{g}+M(n+1)[\sum_{\left|  i\right|
+\left|  j\right|  =n+2,j\neq0}h_{i,j}z_{\ell}p_{i,j-\varepsilon(\ell)}\\
&  +\sum_{\left|  i\right|  =n+2}h_{i,0}\bar{z}_{k}p_{i-\varepsilon
(k),0}]\symbol{94}\\
&  =M(n+1)\hat{g}+M(n+1)W\tilde{h}=M(n+1)(\hat{g}+W\tilde{h})=0\;\;\text{(by
(\ref{zero})).}%
\end{align*}
Thus, $q\in\mathbb{C}_{n+1}^{d}\left[  z,\bar{z}\right]  $ and $q(Z,\bar
{Z})=0$. \ Since $\omega\in\mathcal{V}(M(n+1))$, we must have $q(\omega
,\bar{\omega})=0$. \ Therefore,
\begin{align}
0  &  =g(\omega,\bar{\omega})+\sum_{\left|  i\right|  +\left|  j\right|
=n+2,j\neq0}h_{i,j}\omega_{\ell}p_{i,j-\varepsilon(\ell)}(\omega,\bar{\omega
})\nonumber\\
&  +\sum_{\left|  i\right|  =n+2}h_{i,0}\bar{\omega}_{k}p_{i-\varepsilon
(k),0}(\omega,\bar{\omega}). \label{gw}%
\end{align}
\ Let $r_{i,j}(z,\bar{z}):=\bar{z}^{i}z^{j}-p_{i,j}(z,\bar{z})\;\;(\left|
i\right|  +\left|  j\right|  =n+1).$ \ Clearly each $r_{i,j}\in\mathbb{C}%
_{n+1}^{d}\left[  z,\bar{z}\right]  $ and $r_{i,j}(Z,\bar{Z})=0$ by
(\ref{flatness}), so $r_{i,j}(\omega,\bar{\omega})=0\;\;(\left|  i\right|
+\left|  j\right|  =n+1)$. \ Multiplying $r_{i,j}(\omega,\bar{\omega})=0$ by
either $\omega_{\ell}$ or $\bar{\omega}_{k}$, it follows that
\[
\left\{
\begin{array}
[c]{cc}%
\bar{\omega}^{i}\omega^{j}=(z_{\ell}p_{i,j-\varepsilon(\ell)})(\omega
,\bar{\omega}) & (\left|  i\right|  +\left|  j\right|  =n+2,\;j_{\ell}%
\geq1\text{ for some }\ell=1,...,d\\
\bar{\omega}^{i}=(\bar{z}_{k}p_{i-\varepsilon(k),0})(\omega,\bar{\omega}) &
(\left|  i\right|  =n+2,\;j=0,\;i_{k}\geq1\text{ for some }k=1,...,d
\end{array}
.\ \right.
\]
Now (\ref{gw}) becomes%
\begin{align*}
0  &  =g(\omega,\bar{\omega})+\sum_{\left|  i\right|  +\left|  j\right|
=n+2,j\neq0}h_{i,j}\bar{\omega}^{i}\omega^{j}+\sum_{\left|  i\right|
=n+2}h_{i,0}\bar{\omega}^{i}\\
&  =g(\omega,\bar{\omega})+h(\omega,\bar{\omega})=f(\omega,\bar{\omega}).
\end{align*}
Thus, $f(\omega,\bar{\omega})=0$, as desired.
\end{proof}

\begin{lemma}
\label{W}Assume that $M(n)\equiv M^{d}(n)\geq0$ admits an $r$-atomic
representing measure $\nu$, where $r:=\operatorname*{rank}M(n)$, and let
$\mathcal{V}:=\operatorname*{supp}\nu$. \ If $\mathcal{B}\equiv\{\bar
{Z}^{i_{k}}Z^{j_{k}}\}_{k=1}^{r}$ is a maximal linearly independent subset of
columns of $M(n)$, then $W_{\mathcal{B},\mathcal{V}}$ is invertible (cf.
Theorem \ref{thmB}).
\end{lemma}

\begin{proof}
Let $R_{1},...,R_{r}$ denote the rows of $W_{\mathcal{B},\mathcal{V}}$, and
assume that $W_{\mathcal{B},\mathcal{V}}$ is singular. \ Then there exists
scalars $c_{1},...,c_{r}\in\mathbb{C}$, not all zero, such that $c_{1}%
R_{1}+...+c_{r}R_{r}=0$. \ Let $p(z,\bar{z}):=c_{1}z^{i_{1}}\bar{z}^{j_{1}%
}+...+c_{r}z^{i_{r}}\bar{z}^{j_{r}}$. \ Clearly, $p|_{\operatorname*{supp}\nu
}\equiv0$, so Proposition \ref{MomProp1} implies that $p(Z,\bar{Z})=0$. \ Then
$c_{1}Z^{i_{1}}\bar{Z}^{j_{1}}+...+cZ^{i_{r}}\bar{Z}^{j_{r}}=0$ in
$\mathcal{C}_{M(n)}$, contradicting the fact that $\mathcal{B}$ is linearly independent.
\end{proof}

\begin{proof}
[Proof of Theorem \textup{\ref{thmB}}]Let $r:=\operatorname*{rank}M(n)$; we
first show that $\mathcal{V}$ $\equiv\mathcal{V}(M(n+1))$ satisfies
$\operatorname*{card}\mathcal{V}=r$. \ Theorem \ref{MomThm2} implies that
$M(n+1)$ admits a unique flat extension $M(\infty)$ and that $M(\infty)$
admits an $r$-atomic representing measure $\zeta$. \ Write
$\operatorname*{supp}\zeta\equiv\{\omega_{1},...,\omega_{r}\}$, and define
$p\in\mathbb{C}_{2r}^{d}\left[  z,\bar{z}\right]  $ by $p(z,\bar{z}%
):=\prod_{i=1}^{r}\left\|  z-\omega_{i}\right\|  ^{2}$ (where, for
$z\equiv(z_{1},...,z_{d})$, $\left\|  z\right\|  ^{2}:=\sum_{j=1}^{d}\bar
{z}_{j}z_{j}\in\mathbb{C}_{2}^{d}\left[  z,\bar{z}\right]  )$. \ Clearly,
$\mathcal{Z}(p)=\operatorname*{supp}\zeta$, and since $\zeta$ is a
representing measure for $M(2r)$, Proposition \ref{MomProp1} implies
$p(Z,\bar{Z})=0$ in $\mathcal{C}_{M(2r)}$. \ Thus $\mathcal{V}(M(2r))\subseteq
\mathcal{Z}(p)$ and $\operatorname*{card}\mathcal{V}(M(2r))\leq
\operatorname*{card}\mathcal{Z}(p)=r$. \ To show that $\operatorname*{card}%
\mathcal{V}=r$, we consider two cases. \ If $2r\leq n$, then, since $\zeta$ is
a representing measure for $M(n+1)$, $\operatorname*{supp}\zeta\subseteq
\mathcal{V}(M(n+1))\subseteq\mathcal{V}(M(n))\subseteq\mathcal{V}%
(M(2r))\subseteq\mathcal{Z}(p)=\operatorname*{supp}\zeta$, whence
$\operatorname*{supp}\zeta=\mathcal{V}$ and $\operatorname*{card}%
\mathcal{V}=r$. \ If $2r\geq n+1$, repeated application of Theorem
\ref{propflat} implies that $\mathcal{V}\equiv\mathcal{V}(M(n+1))=\mathcal{V}%
(M(n+2))=...=\mathcal{V}(M(2r))$, and since $\zeta$ is a representing measure
for $M(n+1)$, (\ref{eq.Mom4}) implies
\begin{equation}
r=\operatorname*{rank}M(n+1)\leq\operatorname*{card}\mathcal{V}%
(M(n+1))=...=\operatorname*{card}\mathcal{V}(M(2r)). \label{2a}%
\end{equation}
Now, from above, $\operatorname*{card}\mathcal{V}(M(2r))\leq r$, so (\ref{2a})
implies that $\operatorname*{card}\mathcal{V}=r$ in this case too.

Now let $\nu$ be a representing measure for $M(n+1)$. \ Then
$r=\operatorname*{rank}M(n+1)\leq\operatorname*{card}\operatorname*{supp}%
\nu\leq\operatorname*{card}\mathcal{V}=r$, and since $\operatorname*{supp}%
\nu\subseteq\mathcal{V}$, it follows that $\operatorname*{supp}\nu
=\mathcal{V}$, whence $\nu=\sum_{i=1}^{r}\rho_{i}\delta_{\omega_{i}}$, for
some densities $\rho_{1},...,\rho_{r}$. \ Since $\nu$ is a representing
measure for $M(n)$, $\rho\equiv(\rho_{1},...,\rho_{r})$ satisfies
$W_{\mathcal{B},\mathcal{V}}\rho^{t}=(\gamma_{i_{1},j_{1}},...,\gamma
_{i_{r},j_{r}})^{t}$, and since $W_{\mathcal{B},\mathcal{V}}$ is invertible by
Lemma \ref{W}, $\rho$ is uniquely determined. \ Thus $\nu$ is the unique
representing measure for $M(n+1)$.
\end{proof}

In \cite[Theorem 7.7]{tcmp1} we proved that a finite rank positive infinite
moment matrix $M\equiv M^{d}(\infty)$ has a $\operatorname*{rank}M$-atomic
representing measure; for $d=1$ we established uniqueness in \cite[Theorem
4.7]{tcmp1}. \ We can now establish uniqueness for arbitrary $d$. \ 

\begin{corollary}
\label{newcor}A finite rank positive moment matrix $M\equiv M^{d}(\infty)$ has
a unique representing measure $\nu$, and $\operatorname*{card}%
\operatorname*{supp}\nu=\operatorname*{rank}M$.
\end{corollary}

\begin{proof}
Following \cite[Theorem 7.7]{tcmp1}, let $\zeta$ be a $\operatorname*{rank}%
M$-atomic representing measure for $M$. \ Let $j$ be the smallest integer such
that $\operatorname*{rank}M(j)=\operatorname*{rank}M(j+1)$. \ Theorem
\ref{thmB} implies that $M(j+1)$ has a unique representing measure $\nu$,
whence $\zeta=\nu$ and $\operatorname*{card}\operatorname*{supp}%
\nu=\operatorname*{rank}M.$
\end{proof}

\begin{remark}
The measure $\nu$ in Corollary \ref{newcor} may be computed using Theorem
\ref{thmB}; indeed, $\operatorname*{supp}\nu=\mathcal{V}(M(j+1))$.
\end{remark}

In order to study moment problems on $\mathbb{R}^{N}$, we next introduce real
moment matrices. \ Let $\mathbb{C}^{N}[t]\equiv\mathbb{C}[t_{1},...,t_{N}]$
denote the space of complex polynomials in $N$ real variables, and let
$\mathbb{C}_{s}^{N}[t]$ denote the polynomials of degree at most $s$; \ then
$\dim\mathbb{C}_{s}^{N}[t]=\left(
\begin{array}
[c]{c}%
N+s\\
s
\end{array}
\right)  $. \ For $t\equiv(t_{1},...,t_{N})\in\mathbb{R}^{N}$ and
$i\equiv(i_{1},...,i_{N})\in\mathbb{Z}_{+}^{N}$, we set $t^{i}:=t_{1}^{i_{1}%
}\cdots t_{N}^{i_{N}}$. \ Given a real sequence $\beta\equiv\beta
^{(r)}=\{\beta_{i}\}_{i\in\mathbb{Z}_{+}^{N},\left|  i\right|  \leq r}$, the
truncated moment problem for $\beta$ concerns conditions for the existence of
a positive Borel measure $\mu$ on $\mathbb{R}^{N}$ satisfying
\begin{equation}
\beta_{i}=\int t^{i}\;d\mu(t)\;(\equiv\int t_{1}^{i_{1}}\cdots t_{N}^{i_{N}%
}\;d\mu(t_{1},...,t_{N}))\;\;(\left|  i\right|  \leq r). \label{eq.Mom5}%
\end{equation}
A measure $\mu$ satisfying (\ref{eq.Mom5}) is a \textit{representing measure}
for $\beta$; \ if, in addition, $K\subseteq\mathbb{R}^{N}$ is closed and
$\operatorname*{supp}\mu\subseteq K$, then $\mu$ is a $K$\textit{-representing
measure} for $\beta$.

Let $r=2n$; \ in this case $\beta^{(2n)}$ corresponds to a \textit{real moment
matrix} $\mathcal{M}(n)\equiv\mathcal{M}^{N}(n)(\beta)$, defined as follows.
\ Let $\mathcal{B}\equiv\{t^{i}\}_{i\in\mathbb{Z}_{+}^{N},\left|  i\right|
\leq n}$ denote the basis of monomials in $\mathbb{C}^{N}[t]$, ordered
lexicographically; e.g., for $N=3$, $n=2$, this ordering is $1,t_{1}%
,t_{2},t_{3},t_{1}^{2},t_{1}t_{2},t_{1}t_{3},t_{2}^{2},t_{2}t_{3},t_{3}^{2}$.
\ The size of $\mathcal{M}(n)$ is $\dim\mathbb{C}_{n}^{N}[t]\;(=\left(
\begin{array}
[c]{c}%
N+n\\
n
\end{array}
\right)  )$, with rows and columns indexed as $\{T^{i}\}_{i\in\mathbb{Z}%
_{+}^{N},\left|  i\right|  \leq n}$, following the same lexicographic order as
above. \ The entry of $\mathcal{M}(n)$ in row $T^{i}$, column $T^{j}$ is
$\beta_{i+j}$, $i,j\in\mathbb{Z}_{+}^{N},\left|  i\right|  +\left|  j\right|
\leq2n$. \ Note that for $N=1$, $\mathcal{M}^{N}(n)(\beta)$ is the Hankel
matrix $(\beta_{i+j})$ associated with the classical Hamburger moment problem
($K=\mathbb{R}$) (cf. \cite{Akh}). $\ $

For $p\in\mathbb{C}_{n}^{N}[t]$, $p(t)\equiv\sum_{i\in\mathbb{Z}_{+}%
^{N},\left|  i\right|  \leq n}a_{i}t^{i}$, we let $\tilde{p}:=(a_{i})$ denote
the coefficient vector of $p$ relative to $\mathcal{B}$. \ The Riesz
functional $\Lambda_{\beta}:\mathbb{C}_{2n}^{N}[t]\rightarrow\mathbb{C}$ is
defined by $\Lambda_{\beta}(\sum b_{r}t^{r}):=\sum b_{r}\beta_{r}$. \ Thus,
$\mathcal{M}^{N}(n)(\beta)$ is uniquely determined by
\begin{equation}
\left\langle \mathcal{M}^{N}(n)(\beta)\tilde{f},\tilde{g}\right\rangle
:=\Lambda_{\beta}(f\bar{g})\;\;(f,g\in\mathbb{C}_{n}^{N}[t]). \label{eq.Mom6}%
\end{equation}
If $\beta^{(2n)}$ has a representing measure $\mu$, then $\Lambda_{\beta
}(f\bar{g})=\int f\bar{g}\;d\mu$, so $\mathcal{M}^{N}(n)(\beta)$ is positive
semidefinite. \ 

For $p\equiv\sum_{r\in\mathbb{Z}_{+}^{N},\left|  r\right|  \leq n}a_{r}t^{r}$,
we define an element in $\mathcal{C}_{\mathcal{M}(n)}$ (the column space of
$\mathcal{M}(n)$) by $p(T):=\sum_{r\in\mathbb{Z}_{+}^{N},\left|  r\right|
\leq n}a_{r}T^{r}$. \ Let $\mathcal{V}(\mathcal{M}(n)):=\bigcap
_{\substack{p\in C_{n}^{N}[t]\\p(T)=0}}\mathcal{Z}(p)$ denote the
\textit{variety} of $\mathcal{M}(n)$; we also denote this variety by
$\mathcal{V}(\beta)$. \ Let $J\equiv J(n):=\{j\in\mathbb{Z}_{+}^{N}:\left|
j\right|  \leq n\}$; thus $\operatorname*{card}J(n)=$ size $\mathcal{M}(n)$.
\ Let $s:=$ size $\mathcal{M}(n)-\operatorname*{rank}\mathcal{M}(n)$; the
following result, which proves Proposition \ref{prop1.3}, identifies $s$
polynomials in $\mathbb{R}_{n}^{N}[t]$ whose common zeros comprise
$\mathcal{V}(\mathcal{M}(n))$.

\begin{proposition}
\label{prop1.3new}Let $\mathcal{M}(n)$ be a real moment matrix, with columns
$T^{j}$ indexed by $j\in J$, let $r:=\operatorname*{rank}\mathcal{M}(n)$, and
let $\mathcal{B}\equiv\{T^{i}\}_{i\in I}$ be a maximal linearly independent
set of columns of $\mathcal{M}(n)$, where $I\subseteq J$ satisfies
$\operatorname*{card}I=r$. \ For each index $j\in J\;\backslash\;I$, let
$q_{j}$ denote the unique polynomial in lin.span $\{t^{i}\}_{i\in I}$ such
that $T^{j}=q_{j}(T)$, and let $r_{j}(t):=t^{j}-q_{j}(t)$. \ Then
$\mathcal{V}(\mathcal{M}(n))$ is precisely the set of common zeros of
$\{r_{j}\}_{j\in J\;\backslash\;I}$.
\end{proposition}

\begin{proof}
Clearly $\mathcal{V}\equiv\mathcal{V}(\mathcal{M}(n))\subseteq\bigcap_{j\in
J}\mathcal{Z}(r_{j})$. \ For the reverse inclusion, set $\mathbb{R}_{n}%
^{N}[t]:=\{p\in\mathbb{R}^{N}[t]:\deg p\leq n\}$ and let $\Phi:\mathbb{R}%
_{n}^{N}[t]\rightarrow\mathcal{C}_{\mathcal{M}(n)}$ denote the map $p\mapsto
p(Z,\bar{Z})$. \ $\Phi$ is linear and surjective, so $\dim\ker\Phi
=\dim\mathbb{R}_{n}^{N}[t]-\dim\mathcal{C}_{\mathcal{M}(n)}%
=\operatorname*{card}J-\operatorname*{card}I$. \ Observe now that for $j\in
J\;\backslash\;I$, since $T^{j}=q_{j}(T)$ we have $r_{j}\in\ker\Phi$.
\ Moreover, for $j\in J\;\backslash\;I$, the monomial $t^{j}$ only appears in
$r_{j}$, so it is straightforward to verify that $\{r_{j}\}_{j\in
J\;\backslash\;I}$ is a linearly independent subset of $\mathbb{R}_{n}^{N}%
[t]$. \ It follows at once that $\{r_{j}\}_{j\in J\;\backslash\;I}$ is a basis
for $\ker\Phi$, whence $\bigcap_{j\in J}\mathcal{Z}(r_{j})\subseteq
\bigcap_{p\in\ker\Phi}\mathcal{Z}(p)=\mathcal{V}$.
\end{proof}

\begin{remark}
Proposition \ref{prop1.3new} admits an exact analogue for complex moment matrices.
\end{remark}

We omit the proofs of the following results, which are analogous to the
corresponding proofs for $M^{d}(n)(\gamma)$.

\begin{proposition}
\label{MomProp5}Suppose $\mu$ is a representing measure for $\beta^{(2n)}$.
\ For $p\in\mathbb{C}_{n}^{N}[t]$, $\operatorname*{supp}\mu\subseteq
\mathcal{Z}(p):=\{t\in\mathbb{R}^{N}:p(t)=0\}$ if and only if $p(T)=0$.

\begin{corollary}
\label{MomCor6}If $\beta^{(2n)}$ has a representing measure, then
$\mathcal{M}^{N}(n)(\beta)$ is recursively generated, i.e., if $p,q,pq\in
\mathbb{C}_{n}^{N}[t]$ and $p(T)=0$, then $(pq)(T)=0$.
\end{corollary}
\end{proposition}

\begin{corollary}
\label{MomCor7}If $\mu$ is a representing measure for $\beta^{(2n)}$, then
$\operatorname*{supp}\mu\subseteq\mathcal{V}(\beta)$ and $\operatorname*{card}%
\mathcal{V}(\beta)\geq\operatorname*{card}\operatorname*{supp}\mu
\geq\operatorname*{rank}\mathcal{M}^{N}(n)(\beta)$.
\end{corollary}

We devote the remainder of this section to describing an equivalence between
truncated moment problems on $\mathbb{R}^{2d}$ and $\mathbb{C}^{d}$. \ In the
sequel, $\mathcal{C}^{(n)}$ denotes the ordered basis for $\mathbb{C}_{n}%
^{d}[z,\bar{z}]$ consisting of the monomials, ordered lexicographically by
degree. \ We denote the coefficient vector of $p\in\mathbb{C}_{n}^{d}%
[z,\bar{z}]$ relative to $\mathcal{C}^{(n)}$ by $\hat{p}$; thus $\mathcal{K}%
^{(n)}:=\{\hat{p}:p\in\mathbb{C}_{n}^{d}[z,\bar{z}]\}\cong\mathbb{C}^{\eta
}\cong\mathbb{C}_{n}^{d}[z,\bar{z}]$. \ For $0\leq j\leq n$, let
$\mathcal{K}_{j}$ denote the subspace of $\mathcal{K}^{(n)}$ spanned by
elements $\widehat{\bar{z}^{r}z^{s}}$ with $\left|  r\right|  +\left|
s\right|  =j$; thus $\mathcal{K}^{(n)}=\mathcal{K}^{(n-1)}\bigoplus
\mathcal{K}_{n}\equiv\mathcal{K}_{0}\bigoplus...\bigoplus\mathcal{K}_{n}$, and
$\dim\mathcal{K}_{j}=\binom{j-1+2d}{2d-1}\;(0\leq j\leq n)$.

Next, let $\mathbb{C}_{n}^{2d}[t]\equiv\mathbb{C}_{n}[t_{1},...,t_{2d}]$
denote the vector space over $\mathbb{C}$ of polynomials in real
indeterminates $t_{1},...,t_{2d}$ with total degree at most $n$. \ For
$i\equiv(i_{1},...,i_{2d})\in\mathbb{Z}_{+}^{2d}$, $\left|  i\right|  \leq n$,
let $t^{i}:=t_{1}^{i_{1}}\cdot...\cdot t_{2d}^{i_{2d}}$; thus $q\in
\mathbb{C}_{n}^{2d}[t]$ may be expressed as $q(t)\equiv\sum_{\left|  i\right|
\leq n}b_{i}t^{i}$. \ Note that $\dim\mathbb{C}_{n}^{2d}[t]=\eta(d,n)$. \ In
the sequel, $\mathcal{B}^{(n)}$ denotes the ordered basis for $\mathbb{C}%
_{n}^{2d}[t]$ consisting of the monomials, ordered lexicographically by
degree; for $d=n=2$, this ordering is $1,t_{1},t_{2},t_{3},t_{4},t_{1}%
^{2},t_{1}t_{2},t_{1}t_{3},t_{1}t_{4},t_{2}^{2},t_{2}t_{3},t_{2}t_{4}%
,t_{3}^{2},t_{3}t_{4},t_{4}^{2}$. \ Now we set $x_{i}:=t_{i}\;(1\leq i\leq d)$
and $y_{i}:=t_{i+d}\;(1\leq i\leq d)$, so that $\mathbb{C}_{n}^{2d}%
[t]=\mathbb{C}_{n}^{d}[x,y]:=\mathbb{C}_{n}[x_{1},...,x_{d};y_{1},...,y_{d}]$;
with this notation, for $d=n=2$ the basis $\mathcal{B}^{(2)}$ assumes the form
$1,x_{1},x_{2},y_{1},y_{2},x_{1}^{2},$ $x_{1}x_{2},y_{1}x_{1},y_{2}x_{1}%
,x_{2}^{2},$ $y_{1}x_{2},y_{2}x_{2},y_{1}^{2},y_{1}y_{2},y_{2}^{2}$. \ We
denote the coefficient vector of $q\in\mathbb{C}_{n}^{d}[x,y]$ relative to
$\mathcal{B}^{(n)}$ by $\tilde{q}$; thus $\mathcal{H}^{(n)}:=\{\tilde{q}%
:q\in\mathbb{C}_{n}^{d}[x,y]\}\cong\mathbb{C}^{\eta}\cong\mathbb{C}_{n}%
^{2d}[t]$. \ For $0\leq j\leq n$, let $\mathcal{H}_{j}$ denote the subspace of
$\mathcal{H}^{(n)}$ spanned by elements $\widetilde{y^{r}x^{s}}$ with $\left|
r\right|  +\left|  s\right|  =j$; thus $\mathcal{H}^{(n)}=\mathcal{H}%
^{(n-1)}\bigoplus\mathcal{H}_{n}\equiv\mathcal{H}_{0}\bigoplus...\bigoplus
\mathcal{H}_{n}$, and $\dim\mathcal{H}_{j}=\binom{j-1+2d}{2d-1}\;(0\leq j\leq
n)$.

For $0\leq j\leq n$, we define a linear map $L_{j}:\mathcal{K}_{j}%
\rightarrow\mathcal{H}_{j}$ by $L_{j}(\widehat{\bar{z}^{k}z^{\ell}%
}):=[(x-iy)^{k}(x+iy)^{\ell}]^{\symbol{126}}\;(\left|  k\right|  +\left|
\ell\right|  =j)$. \ Since $(x-iy)^{k}(x+iy)^{\ell}\equiv(x_{1}-iy_{1}%
)^{k_{1}}\cdot...\cdot(x_{d}-iy_{d})^{k_{d}}(x_{1}+iy_{1})^{\ell_{1}}%
\cdot...\cdot(x_{d}+iy_{d})^{\ell_{d}}$, the Binomial Theorem shows that
$L_{j}(\widehat{\bar{z}^{k}z^{\ell}})$ is indeed an element of $\mathcal{H}%
_{j}$. \ We now define $L\equiv L^{(n)}:\mathcal{K}^{(n)}\rightarrow
\mathcal{H}^{(n)}$ by $L:=\bigoplus_{k=0}^{n}L_{k}\;(=L^{(n-1)}\bigoplus
L_{n}).$ For $d=n=2$, we have
\[
L_{0}=(1),\;L_{1}=\left(
\begin{array}
[c]{cccc}%
1 & 0 & 1 & 0\\
0 & 1 & 0 & 1\\
i & 0 & -i & 0\\
0 & i & 0 & -i
\end{array}
\right)  ,
\]%
\[
L_{2}=\left(
\begin{array}
[c]{cccccccccc}%
1 & 0 & 1 & 0 & 0 & 0 & 0 & 1 & 0 & 0\\
0 & 1 & 0 & 1 & 0 & 1 & 0 & 0 & 1 & 0\\
2i & 0 & 0 & 0 & 0 & 0 & 0 & -2i & 0 & 0\\
0 & i & 0 & -i & 0 & i & 0 & 0 & -i & 0\\
0 & 0 & 0 & 0 & 1 & 0 & 1 & 0 & 0 & 1\\
0 & i & 0 & i & 0 & -i & 0 & 0 & -i & 0\\
0 & 0 & 0 & 0 & 2i & 0 & 0 & 0 & 0 & -2i\\
-1 & 0 & 1 & 0 & 0 & 0 & 0 & -1 & 0 & 0\\
0 & -1 & 0 & 1 & 0 & 1 & 0 & 0 & -1 & 0\\
0 & 0 & 0 & 0 & -1 & 0 & 1 & 0 & 0 & -1
\end{array}
\right)
\]
and $L^{(2)}=L_{0}\bigoplus L_{1}\bigoplus L_{2}=L^{(1)}\bigoplus L_{2}$.
$\ $To clarify the properties of $L$ we introduce the map $\psi:\mathbb{R}%
^{d}\times\mathbb{R}^{d}\rightarrow\mathbb{C}^{d}\times\mathbb{C}^{d}$ defined
by $\psi(x,y):=(z,\bar{z})$, where $z\equiv x+iy$, $\bar{z}\equiv
x-iy\in\mathbb{C}^{d}$. \ Clearly $\psi$ is injective, and we let
\newline $\tau:$ Ran$\;\psi\rightarrow\mathbb{R}^{d}\times\mathbb{R}^{d}$
denote the inverse map, $\tau(z,\bar{z}):=(\frac{z+\bar{z}}{2},\frac{z-\bar
{z}}{2i})$.

\begin{lemma}
\label{lem10}(i) $L\hat{p}:=\widetilde{p\circ\psi}$ ($p\in\mathbb{C}_{n}%
^{d}[z,\bar{z}]$).\newline (ii) $L$ is invertible, with $L^{-1}(\tilde
{q})=(q\circ\tau)\symbol{94}$.
\end{lemma}

\begin{proof}
(i) For $p\in\mathbb{C}_{n}^{d}[z,\bar{z}]$, write $p(z,\bar{z})\equiv
\sum_{\left|  k\right|  +\left|  \ell\right|  \leq n}a_{k\ell}\bar{z}%
^{k}z^{\ell}$. \ Then
\begin{align*}
L(\hat{p})  &  =\sum_{\left|  k\right|  +\left|  \ell\right|  \leq n}a_{k\ell
}L(\widehat{\bar{z}^{k}z^{\ell}})=\sum_{\left|  k\right|  +\left|
\ell\right|  \leq n}a_{k\ell}[(x-iy)^{k}(x+iy)^{\ell}]^{\symbol{126}}\\
&  =\sum_{\left|  k\right|  +\left|  \ell\right|  \leq n}a_{k\ell}[\bar{z}%
^{k}z^{\ell}\circ\psi]^{\symbol{126}}=[(\sum_{\left|  k\right|  +\left|
\ell\right|  \leq n}a_{k\ell}\bar{z}^{k}z^{\ell})\circ\psi]^{\symbol{126}%
}=\widetilde{p\circ\psi}.
\end{align*}
\newline (ii) A calculation shows that $R_{j}:=L_{j}^{-1}:\mathcal{H}%
_{j}\rightarrow\mathcal{K}_{j}$ is given by $R_{j}(\widetilde{y^{r}x^{s}%
}):=[(\frac{z-\bar{z}}{2i})^{r}(\frac{z+\bar{z}}{2})^{s}]^{\symbol{94}}$
$\equiv\lbrack(\frac{z_{1}-\bar{z}_{1}}{2i})^{r_{1}}\cdot...\cdot
(\frac{z_{d}-\bar{z}_{d}}{2i})^{r_{d}}(\frac{z_{1}+\bar{z}_{1}}{2})^{s_{1}%
}\cdot...\cdot(\frac{z_{d}+\bar{z}_{d}}{2})^{s_{d}}]^{\symbol{94}}$; thus
$L^{-1}:\mathcal{H}^{(n)}\rightarrow\mathcal{K}^{(n)}$ satisfies
$L^{-1}(\tilde{q})=(q\circ\tau)\symbol{94}$.
\end{proof}

Our next goal is to associate to a complex sequence $\gamma\equiv\gamma
^{(2n)}=\{\gamma_{rs}\}_{r,s\in\mathbb{Z}_{+}^{d},\left|  r\right|  +\left|
s\right|  \leq2n}$, with $\gamma_{00}>0$ and $\gamma_{rs}=\bar{\gamma}_{sr}$,
an ``equivalent'' real sequence $\beta\equiv\beta^{(2n)}=\{\beta_{j}%
\}_{j\in\mathbb{Z}_{+}^{2d},\left|  j\right|  \leq2n}$, with $\beta_{0}%
=\gamma_{00}$. \ We require the following lemma.

\begin{lemma}
\label{real}Let $p\equiv\sum a_{rs}\bar{z}^{r}z^{s}\in\mathbb{C}_{2n}%
^{d}[z,\bar{z}]$ and assume that $p$ is real-valued. \ Then $\Lambda_{\gamma
}(p)$ is real.
\end{lemma}

\begin{proof}
Recall that $\Lambda_{\gamma}(p)=\sum a_{rs}\gamma_{rs}$. \ Then
$\overline{\Lambda_{\gamma}(p)}=\sum\bar{a}_{rs}\bar{\gamma}_{rs}%
=\sum\bar{a}_{rs}\gamma_{sr}=\Lambda_{\gamma}(\bar{p})=\Lambda_{\gamma}(p)$,
so $\Lambda_{\gamma}(p)$ is real.
\end{proof}

For $j\in\mathbb{Z}_{+}^{2d}$, $\left|  j\right|  \leq2n$, set $\pi
_{x}(j):=(j_{1},...,j_{d})$ and $\pi_{y}(j):=(j_{d+1},...,j_{2d})$. \ For
$\gamma$ as above, we now set $\beta_{j}:=\Lambda_{\gamma}(y^{\pi_{y}%
(j)}x^{\pi_{x}(j)})$, where, for $z\in\mathbb{C}^{d}$, $x:=\frac{z+\bar{z}}%
{2}$ and $y:=\frac{z-\bar{z}}{2i}$. \ Since the operand of $\Lambda_{\gamma}$
is real-valued (as an element of $\mathbb{C}_{2n}^{d}[z,\bar{z}]$), Lemma
\ref{real} implies $\beta_{j}\in\mathbb{R}$. \ We now set $\mathcal{R}%
(\gamma):=\beta$; note that
\begin{equation}
\beta_{j}=\Lambda_{\beta}(t^{j})=\Lambda_{\beta}(y^{\pi_{y}(j)}x^{\pi_{x}%
(j)})=\Lambda_{\gamma}((\frac{z-\bar{z}}{2i})^{\pi_{y}(j)}(\frac{z+\bar{z}}%
{2})^{\pi_{x}(j)}). \label{new12}%
\end{equation}

\begin{proposition}
\label{prop12}$\mathcal{M}(n)(\mathcal{R}(\gamma))=L^{\ast-1}M(n)(\gamma
)L^{-1}$.
\end{proposition}

\begin{proof}
It suffices to show that for $k,\ell,r,s\in\mathbb{Z}_{+}^{d}$, with $\left|
k\right|  +\left|  \ell\right|  ,\left|  r\right|  +\left|  s\right|  \leq n$,
and for $\beta=\mathcal{R}(\gamma)$, we have $\left\langle \mathcal{M}%
(n)(\beta)\widetilde{y^{k}x^{\ell}},\widetilde{y^{r}x^{s}}\right\rangle
=\left\langle L^{\ast-1}M(n)(\gamma)L^{-1}\widetilde{y^{k}x^{\ell}}%
,\widetilde{y^{r}x^{s}}\right\rangle $. \ Now,
\begin{align}
\left\langle L^{\ast-1}M(n)(\gamma)L^{-1}\widetilde{y^{k}x^{\ell}}%
,\widetilde{y^{r}x^{s}}\right\rangle  &  =\left\langle M(n)(\gamma
)L^{-1}\widetilde{y^{k}x^{\ell}},L^{-1}\widetilde{y^{r}x^{s}}\right\rangle
\nonumber\\
&  =\left\langle M(n)(\gamma)[(\frac{z-\bar{z}}{2i})^{k}(\frac{z+\bar{z}}%
{2})^{\ell}]\symbol{94},[(\frac{z-\bar{z}}{2i})^{r}(\frac{z+\bar{z}}{2}%
)^{s}]\symbol{94}\right\rangle \;\;\nonumber\\
&  \text{(by Lemma \ref{lem10})}\nonumber\\
&  =\Lambda_{\gamma}((\frac{z-\bar{z}}{2i})^{k+r}(\frac{z+\bar{z}}{2}%
)^{\ell+s}). \label{new13}%
\end{align}
Choosing $j\in\mathbb{Z}_{+}^{2d}$ so that $\pi_{x}(j)=\ell+s$ and $\pi
_{y}(j)=k+r$, we have $\left|  j\right|  =(\left|  k\right|  +\left|
\ell\right|  )+(\left|  r\right|  +\left|  s\right|  )\leq2n$, so
(\ref{new12}) shows that the expression in (\ref{new13}) is equal to
$\Lambda_{\beta}(y^{k+r}x^{\ell+s})=\left\langle \mathcal{M}(n)(\beta
)\widetilde{y^{k}x^{\ell}},\widetilde{y^{r}x^{s}}\right\rangle $, as desired.
\end{proof}

Next, we define an inverse to $\mathcal{R}$. $\ $Given a real sequence
$\beta\equiv\beta^{(2n)}=\{\beta_{j}\}_{j\in\mathbb{Z}_{+}^{2d},\left|
j\right|  \leq2n}$, with $\beta_{0}>0$, we will associate to $\beta$ a complex
sequence $\gamma\equiv\gamma^{(2n)}$. $\ $For $k,\ell\in\mathbb{Z}_{+}^{d}$,
$\left|  k\right|  +\left|  \ell\right|  \leq2n$, let
\begin{align*}
\gamma_{k\ell}  &  :=\Lambda_{\beta}((x-iy)^{k}(x+iy)^{\ell})\\
&  =\Lambda_{\beta}((t_{1}-it_{d+1})^{k_{1}}\cdot...\cdot(t_{d}-it_{2d}%
)^{k_{d}}(t_{1}+it_{d+1})^{\ell_{1}}\cdot...\cdot(t_{d}+it_{2d})^{\ell_{d}}).
\end{align*}
Clearly, $\gamma_{00}=\Lambda_{\beta}(1)=\beta_{0}>0$, and $\gamma_{\ell
k}=\bar{\gamma}_{k\ell}$. \ We set $\mathcal{S}(\beta):=\gamma$; we omit the
proof of the following result, which is dual to that in Proposition
\ref{prop12}.

\begin{proposition}
\label{prop13}$M(n)(\mathcal{S}(\beta))=L^{\ast}\mathcal{M}(n)(\beta)L$.
\end{proposition}

Taken together, Propositions \ref{prop12} and \ref{prop13} show that
$(\mathcal{R}\circ\mathcal{S})(\beta)=\beta$ and $(\mathcal{S}\circ
\mathcal{R})(\gamma)=\gamma$. \ We are now in position to formulate the
equivalence between the real and complex truncated moment problems, as
expressed in the following two results.

\begin{proposition}
\label{PropositionEquiv}Given $\gamma\equiv\gamma^{(2n)}$, let $\beta
\equiv\beta^{(2n)}:=\mathcal{R}(\gamma)$. \ \newline (i) $\mathcal{M}%
(n)(\beta)=L^{\ast-1}M(n)(\gamma)L^{-1}$.\newline (ii) $\mathcal{M}%
(n)(\beta)\geq0\Leftrightarrow M(n)(\gamma)\geq0$.\newline (iii)
$\operatorname*{rank}\mathcal{M}(n)(\beta)=\operatorname*{rank}M(n)(\gamma
)$.\newline (iv) $\mathcal{M}(n)(\beta)$ is positive and admits a flat
extension $\mathcal{M}(n+1)$ if and only if $M(n)(\gamma)$ is positive and
admits a flat extension $M(n+1)$.\newline (v) For $q\in\mathbb{C}_{n}[x,y]$,
$q(X,Y)=L^{\ast-1}((q\circ\tau)(Z,\bar{Z}))$.\newline (vi) For $q\in
\mathbb{C}_{n}[x,y]$, $\Lambda_{\beta}(q)=\Lambda_{\gamma}(q\circ\tau
)$\newline (vii) If $\nu$ is a representing measure for $\gamma$, then
$\mu:=\nu\circ\psi$ is a representing measure for $\beta$, of the same measure
class and cardinality of support; moreover, $\operatorname*{supp}\mu
=\tau(\operatorname*{supp}\nu)$.
\end{proposition}

\begin{proof}
(i) This is Proposition \ref{prop12}.

(ii) This follows from (i) and the invertibility of $L$ (Lemma \ref{lem10}).

(iii) This also follows from (i) and the invertibility of $L$.

(iv) Suppose $M(n)(\gamma)$ is positive and admits a flat extension
\[
M(n+1)(\tilde{\gamma})\equiv\left(
\begin{array}
[c]{cc}%
M(n)(\gamma) & B(n+1)\\
B(n+1)^{\ast} & C(n+1)
\end{array}
\right)  .\
\]
Proposition \ref{prop12} (using $n+1$) implies that $\mathcal{M}%
:=(L^{(n+1)\ast})^{-1}M(n+1)(\tilde{\gamma})(L^{(n+1)})^{-1}$ is of the form
$M(n+1)(\mathcal{R}(\tilde{\gamma}))$, while (i) and the direct sum structure
of $(L^{(n+1)})^{-1}$ show that
\[
\mathcal{M}=\left(
\begin{array}
[c]{cc}%
(L^{(n)\ast})^{-1}M(n)(\gamma)(L^{(n)})^{-1} & \ast\\
\ast & \ast
\end{array}
\right)  =\left(
\begin{array}
[c]{cc}%
\mathcal{M}(n)(\mathcal{R}(\gamma)) & \ast\\
\ast & \ast
\end{array}
\right)  .\
\]
Since $\operatorname*{rank}\mathcal{M}=\operatorname*{rank}M(n+1)(\tilde
{\gamma})=\operatorname*{rank}M(n)(\gamma)=\operatorname*{rank}\mathcal{M}%
(n)(\beta)$, it follows that $\mathcal{M}$ is a flat extension of
$\mathcal{M}(n)(\beta)\;(\geq0)$. \ The converse is proved similarly, using
Proposition \ref{prop13}; we omit the details.%
\begin{align*}
\text{(v)\ \ \ \ \ \ \ \ \ \ \ \ \ \ \ \ \ \ \ \ \ \ \ \ }q(X,Y)  &
\equiv\mathcal{M}(n)(\beta)\tilde{q}=L^{\ast}{}^{-1}M(n)(\gamma)L^{-1}%
\tilde{q}\;\;\text{(by (i))}\\
&  =L^{\ast}{}^{-1}M(n)(\gamma)\widehat{q\circ\tau}\;\;\text{(by Lemma
\ref{lem10})}\\
&  =L^{\ast}{}^{-1}(q\circ\tau)(Z,\bar{Z}).
\end{align*}
(vi) Straightforward from (\ref{new12}).

(vii) For $j\in\mathbb{Z}_{+}^{2d}$, $\left|  j\right|  \leq2n$,
\begin{align*}
\int t^{j}\;d\mu &  =\int y^{\pi_{y}(j)}x^{\pi_{x}(j)}\;d(\nu\circ\psi)(x,y)\\
&  =\int(\frac{z-\bar{z}}{2i})^{\pi_{y}(j)}(\frac{z+\bar{z}}{2})^{\pi_{x}%
(j)}\;d\nu(z,\bar{z})\\
&  =\Lambda_{\gamma}((\frac{z-\bar{z}}{2i})^{\pi_{y}(j)}(\frac{z+\bar{z}}%
{2})^{\pi_{x}(j)})\\
&  =\beta_{j}\;\;\text{(by (\ref{new12}));}%
\end{align*}
thus, $\mu$ is a representing measure for $\beta$, and the other properties of
$\mu$ are clear.
\end{proof}

We omit the proof of the following result, which is dual to Proposition
\ref{PropositionEquiv}.

\begin{proposition}
\label{prop15}Given $\beta\equiv\beta^{(2n)}$, let $\gamma\equiv\gamma
^{(2n)}:=\mathcal{S}(\beta)$. \ \newline (i) $M(n)(\gamma)=L^{\ast}%
\mathcal{M}(n)(\beta)L$.\newline (ii) $M(n)(\gamma)\geq0\Leftrightarrow
\mathcal{M}(n)(\beta)\geq0$.\newline (iii) $\operatorname*{rank}%
M(n)(\gamma)=\operatorname*{rank}\mathcal{M}(n)(\beta)$.\newline (iv)
$M(n)(\gamma)$ is positive and admits a flat extension $M(n+1)$ if and only if
$\mathcal{M}(n)(\beta)$ is positive and admits a flat extension $\mathcal{M}%
(n+1)$.\newline (v) For $p\in\mathbb{C}_{n}[z,\bar{z}]$, $p(Z,\bar{Z}%
)=L^{\ast}((p\circ\psi)(X,Y))$.\newline (vi) For $p\in\mathbb{C}_{n}[z,\bar
{z}]$, $\Lambda_{\gamma}(p)=\Lambda_{\beta}(p\circ\psi)$.\newline (vii) If
$\mu$ is a representing measure for $\beta$, then $\nu:=\mu\circ\tau$ is a
representing measure for $\gamma$, of the same measure class and cardinality
of support; moreover, $\operatorname*{supp}\nu=\psi(\operatorname*{supp}\mu)$.
\end{proposition}

Throughout the sequel, whenever we have equivalent sequences $\gamma$ and
$\beta$ (as described by the preceding results), the context always indicates
whether we have $\beta=\mathcal{R}(\gamma)$ or $\gamma=\mathcal{S}(\beta)$, so
we do not explicitly refer to $\mathcal{R}$ or $\mathcal{S}$.

We next present an analogue of Theorem \ref{MomThm2} for truncated moment
problems on $\mathbb{R}^{N}$.

\begin{theorem}
\label{newthm}Let $\beta\equiv\beta^{(2n)}$ and let $r:=\operatorname*{rank}%
\mathcal{M}^{N}(n)(\beta)$. \ If $\mu$ is an $r$--atomic representing measure
for $\beta$, then $\mathcal{M}^{N}(n+1)[\mu]$ is a flat (positive) extension
of $\mathcal{M}(n)\equiv\mathcal{M}^{N}(n)(\beta)$. \ Conversely, if
$\mathcal{M}(n)$ is positive semidefinite and admits a flat extension
$\mathcal{M}(n+1)\equiv\mathcal{M}^{N}(n+1)(\tilde{\beta})$, then
$\mathcal{M}(n+1)$ admits unique flat positive moment matrix extensions
$\mathcal{M}^{N}(n+2)(\tilde{\beta})$, $\mathcal{M}^{N}(n+3)(\tilde{\beta}%
)$,..., and there exists an $r$-atomic representing measure for $\mathcal{M}%
^{N}(\infty)(\tilde{\beta})$ (i.e., a representing measure for $\tilde{\beta
}^{(\infty)}$).
\end{theorem}

\begin{proof}
Suppose $\mu$ is an $r$-atomic representing measure for $\beta$, i.e.,
$\mathcal{M}^{N}(n)(\beta)=\mathcal{M}^{N}(n)[\mu]$. \ Since $\mu$ is also a
representing measure for $\mathcal{M}^{N}(n+1)[\mu]$, Corollary \ref{MomCor7}
implies that $r=\operatorname*{card}\operatorname*{supp}\mu\geq
\operatorname*{rank}\mathcal{M}^{N}(n+1)[\mu]\geq\operatorname*{rank}%
\mathcal{M}^{N}(n)[\mu]=r$, so $\mathcal{M}^{N}(n+1)[\mu]$ is a flat
(positive) extension of $\mathcal{M}(n)$.

For the converse, we assume that $\mathcal{M}^{N}(n)(\beta)$ is positive and
admits a flat extension $\mathcal{M}^{N}(n+1)(\tilde{\beta})$. \ We consider
first the case when $N$ is even, say $N=2d$. \ In this case, let $\gamma
\equiv\gamma^{(2n)}=\mathcal{S}(\beta)$. \ Proposition \ref{prop15} implies
that $M^{d}(n)(\gamma)$ is positive and admits a flat extension $M^{d}%
(n+1)(\tilde{\gamma})$. Theorem \ref{MomThm2} now implies that $M^{d}%
(n+1)(\tilde{\gamma})$ admits unique successive flat (positive) extensions
$M^{d}(n+2)(\tilde{\gamma}),M^{d}(n+3)(\tilde{\gamma}),...,$ and that
$\tilde{\gamma}^{(\infty)}$ admits an $r$-atomic representing measure $\nu$.
\ Proposition \ref{PropositionEquiv} (and the direct sum structure of
$L^{(n+j)}\;(j\geq0)$) now imply that $\mathcal{M}^{2d}(n+1)(\tilde{\beta})$
admits unique successive flat extensions $\{\mathcal{M}^{2d}(n+j)(\tilde
{\beta})\}_{j\geq2}$, defined by $\mathcal{M}^{2d}(n+j)(\tilde{\beta
}):=(L^{(n+j)\ast})^{-1}M^{d}(n+j)(\tilde{\gamma})(L^{(n+j)})^{-1}$.
\ Proposition \ref{PropositionEquiv} further implies that $\nu$ corresponds to
an $r$-atomic representing measure $\mu$ for $\tilde{\beta}^{(\infty)}$.

We now consider the case $N=2d-1$. \ For $x\in\mathbb{R}^{2d-1}$,
$t\in\mathbb{R}$, $i\in\mathbb{Z}_{+}^{2d-1}$, $j\in\mathbb{Z}_{+}$, we set
$\breve{x}:=(x,t)\in\mathbb{R}^{2d}$ and $\breve{\imath}:=(i,j)\in
\mathbb{Z}_{+}^{2d}$, so that $\breve{x}^{\breve{\imath}}=x^{i}t^{j}$.
\ Corresponding to $\beta\equiv\beta^{(2n)}=\{\beta_{i}\}_{i\in\mathbb{Z}%
_{+}^{2d-1}},\left|  i\right|  \leq2n$, we define a sequence $\breve{\beta
}\equiv\breve{\beta}^{(2n)}=\{\breve{\beta}_{\breve{\imath}}\}_{\breve{\imath
}\in\mathbb{Z}_{+}^{2d},\left|  \breve{\imath}\right|  \leq2n}$ as follows:
\begin{equation}
\breve{\beta}_{\breve{\imath}}:=\left\{
\begin{array}
[c]{cc}%
\beta_{i} & \text{if }j=0\\
0 & \text{if }j>0
\end{array}
\right.  .\ \label{betadef}%
\end{equation}
Corresponding to $\mathcal{M}\equiv\mathcal{M}^{2d-1}(n)(\beta)$ we define the
moment matrix $\mathcal{\breve{M}}\equiv\mathcal{M}^{2d}(n)(\breve{\beta})$.
\ Since $\mathcal{\breve{M}}$ is unitarily equivalent to a matrix of the form
$\mathcal{M}\bigoplus0$, we have $\operatorname*{rank}\mathcal{\breve{M}%
=}\operatorname*{rank}\mathcal{M}$, and $\mathcal{\breve{M}}\geq0$ if and only
if $\mathcal{M}\geq0$. \ Suppose $\mathcal{M}\equiv\mathcal{M}^{2d-1}%
(n)(\beta)\geq0$ and suppose $\mathcal{M}(n+1)\equiv\mathcal{M}^{2d-1}%
(n+1)(\tilde{\beta})$ is a flat extension of $\mathcal{M}$. $\ $We claim that
$\mathcal{\breve{M}}(n+1)\equiv\lbrack M(n+1)]\symbol{94}$ is a flat extension
of $\mathcal{\breve{M}}$. \ Since $\mathcal{M}(n+1)\geq0$, then
$\mathcal{\breve{M}}(n+1)\geq0$, and $\operatorname*{rank}\mathcal{\breve{M}%
}(n+1)=\operatorname*{rank}\mathcal{M}(n+1)=\operatorname*{rank}%
\mathcal{M}(n)=\operatorname*{rank}\mathcal{\breve{M}}(n)$. \ Let us denote
$\mathcal{\breve{M}}(n+1)$ as $\mathcal{M}^{2d}(n+1)(\lambda)$, for some
sequence $\lambda$. \ To show that $\mathcal{\breve{M}}(n+1)$ is an extension
of $\mathcal{\breve{M}}(n)$, it suffices to show that if $\breve{\imath}$
satisfies $\left|  \breve{\imath}\right|  \leq2n$, then $\lambda
_{\breve{\imath}}=\breve{\beta}_{\breve{\imath}}$. \ Indeed, if $\breve
{\imath}=(i,j)$ and $j=0$, then $\lambda_{\breve{\imath}}=\tilde{\beta}%
_{i}=\beta_{i}=\breve{\beta}_{\breve{\imath}}$, while if $j>0$, then
$\lambda_{\breve{\imath}}=0=\breve{\beta}_{\breve{\imath}}$. \ Thus
$\mathcal{\breve{M}}(n+1)$ is a flat (positive) extension of $\mathcal{\breve
{M}}(n)$.

Since $\mathcal{\breve{M}}(n+1)=\mathcal{M}^{2d}(n+1)(\lambda)$, the ``even''
case (above) implies that $\mathcal{\breve{M}}(n+1)$ has unique successive
flat moment matrix extensions $\mathcal{M}^{2d}(n+j)(\tilde{\lambda}%
)\;(j\geq2)$, and that $\tilde{\lambda}^{(\infty)}$ admits a
$\operatorname*{rank}\mathcal{\breve{M}}$-atomic representing measure $\nu$.
\ For $j\geq2$ and $i\in\mathbb{Z}_{+}^{2d-1}$ with $\left|  i\right|
\leq2(n+j)$, we set $\tilde{\beta}_{i}:=\tilde{\lambda}_{(i,0)}$. \ Then
$\mathcal{M}^{2d-1}(n+2)(\tilde{\beta}),\mathcal{M}^{2d-1}(n+3)(\tilde{\beta
}),...$, define the unique successive flat moment matrix extensions of
$\mathcal{M}^{2d-1}(n+1)(\tilde{\beta})$ (indeed, $[\mathcal{M}^{2d-1}%
(n+1)(\tilde{\beta})]^{\smallsmile}=M^{d}(n+j)(\tilde{\lambda})\;(j\geq1)$).
\ Finally, if $\nu\equiv\sum_{s=1}^{r}\rho_{s}\delta_{(x_{s},t_{s})}\;$(with
$x_{s}\in\mathbb{R}^{2d-1}$, $t_{s}\in\mathbb{R}$, $\rho_{s}>0$), then
$\mu:=\sum_{s=1}^{r}\rho_{s}\delta_{x_{s}}$ is an $r$-atomic representing
measure for $\tilde{\beta}^{(\infty)}$.
\end{proof}

\begin{remark}
\label{Exi.Rem4}We note the following for future reference. \ In
$\mathcal{\breve{M}}(n+1)\equiv\mathcal{M}^{2d}(n+1)(\lambda)$, since
$\lambda_{\breve{\imath}}=0$ whenever $\left|  \breve{\imath}\right|
\leq2(n+1)$ and $j>0$, each column that is indexed by a multiple of $t$ is
identically $0$. \ Further, since $\tilde{\lambda}^{(\infty)}$ has a
representing measure, each of the successive flat extensions $\mathcal{M}%
^{2d}(n+j)(\tilde{\lambda})\;(j\geq2)$ is recursively generated; hence, in
$\mathcal{M}^{2d}(n+j)(\tilde{\lambda})$, each column indexed by a multiple of
$t$ is identically $0$, whence $\tilde{\lambda}_{(i,j)}=0$ whenever $j>0$.
\end{remark}

We can now give a proof of Theorem \ref{thm12}, which we restate here for the
reader's convenience.

\begin{theorem}
\label{thmA}If $\mathcal{M}(n)\equiv\mathcal{M}^{N}(n)(\beta)\geq0$ admits a
flat extension $\mathcal{M}(n+1)$, then $\mathcal{V}:=\mathcal{V}%
(\mathcal{M}(n+1))$ satisfies $\operatorname*{card}\mathcal{V}=r\;(\equiv
\operatorname*{rank}\mathcal{M}(n))$, and $\mathcal{V}\equiv\{t_{j}%
\}_{j=1}^{r}\subseteq\mathbb{R}^{N}$ forms the support of the unique
representing measure $\mu$ for $\mathcal{M}(n+1)$. \ If $\mathcal{B}%
\equiv\{T^{i_{k}}\}_{k=1}^{r}$ is a maximal linearly independent subset of
columns of $\mathcal{M}(n)$, then $W_{\mathcal{B},\mathcal{V}}$ is invertible,
and $\mu=\sum_{i=1}^{r}\rho_{j}\delta_{t_{j}}$, where $\rho\equiv(\rho
_{1},...,\rho_{r})$ is uniquely determined by $\rho^{t}=W_{\mathcal{B}%
,\mathcal{V}}^{-1}(\beta_{i_{1}},...,\beta_{i_{r}})^{t}$.
\end{theorem}

\begin{proof}
We first consider the support of a representing measure $\mu$ for
$\mathcal{M}(n+1)$ (cf. Theorem \ref{newthm}). \ For $N=2d$, let $\gamma$ be
the equivalent complex sequence associated to $\beta$ via Proposition
\ref{prop15}; Propositions \ref{PropositionEquiv}(v) and \ref{prop15}(v) imply
that $\mathcal{V}(\mathcal{M}(n+1)(\beta))$ and $\mathcal{V}(M(n+1)(\gamma))$
are identical when regarded as subsets of $\mathbb{R}^{2d}$. \ The conclusion
that $\operatorname*{card}\mathcal{V}=r$ and $\operatorname*{supp}%
\mu=\mathcal{V}$ thus follows by a straightforward application of Theorem
\ref{thmB} and Propositions \ref{PropositionEquiv} and \ref{prop15}. \ For
$N=2d-1$, one needs to argue as in the proof of Theorem \ref{newthm}, to
convert the initial moment problem for $\beta$ into an equivalent one for
$\hat{\beta}$ in $\mathbb{R}^{2d}$ (using (\ref{betadef})), and to then appeal
to the result for $N=2d$. \ We omit the details of this argument, except to
note that in the notation of the proof of Theorem \ref{newthm}, $\mathcal{V}%
(\mathcal{M}(n+1)(\beta))\times\{0\}=\mathcal{V}(\mathcal{M}(n+1)(\lambda))$.
\ As for the uniqueness of $\mu$ and the calculation of the densities using
$W_{\mathcal{B},\mathcal{V}}$, the proof is very similar to the argument
establishing the uniqueness of $\nu$ in Theorem \ref{thmB}; for this we use an
analogue of Lemma \ref{W} for the invertibility of $W_{\mathcal{B}%
,\mathcal{V}}$ in the case of real moment matrices.
\end{proof}

\begin{remark}
Theorem \ref{propflat} and Corollary \ref{newcor} admit exact analogues for
real moment matrices.
\end{remark}

\section{\label{Loc}Localizing matrices}

Let $1\leq k\leq n$ and let $p\equiv p\left(  z,\bar{z}\right)  \in
\mathbb{C}^{d}\left[  z,\bar{z}\right]  $, $\deg p=2k$ or $2k-1$. We next
define the \textit{localizing matrix} $M_{p}^{d}\left(  n\right)  \equiv
M_{p}^{d}\left(  n\right)  \left(  \gamma\right)  $ whose positivity is
directly related to the existence of a representing measure for $\gamma
\equiv\gamma^{\left(  2n\right)  }$ with support in $\mathcal{K}_{p}%
\equiv\left\{  z\in\mathbb{C}^{d}:p\left(  z,\bar{z}\right)  \geq0\right\}  $.
Note that $\dim\mathbb{C}_{n-k}^{d}\left[  z,\bar{z}\right]  =\eta\equiv
\eta\left(  d,n-k\right)  =\binom{n-k+2d}{2d}$; thus $\mathbb{C}^{\eta
}=\left\{  \hat{f}:f\in\mathbb{C}_{n-k}^{d}\left[  z,\bar{z}\right]  \right\}
$. We define the $\eta\times\eta$ matrix $M_{p}^{d}\left(  n\right)  $ by%
\begin{equation}
\left\langle M_{p}^{d}\left(  n\right)  \hat{f},\hat{g}\right\rangle
=\Lambda_{\gamma}\left(  pf\bar{g}\right)  \qquad(f,g\in\mathbb{C}_{n-k}%
^{d}\left[  z,\bar{z}\right]  ). \label{eqLoc.1}%
\end{equation}
If $\gamma$ has a representing measure $\nu$ supported in $\mathcal{K}_{p}$,
then $\left\langle M_{p}\left(  n\right)  \hat{f},\hat{f}\right\rangle
=\Lambda_{\gamma}\left(  p\left|  f\right|  ^{2}\right)  =\int p\left|
f\right|  ^{2}\,d\nu\geq0$, whence $M_{p}^{d}\left(  n\right)  \geq0$. Note
also the following consequences of (\ref{eqLoc.1}):%
\begin{equation}
M_{p}^{d}\left(  n\right)  ^{\ast}=M_{\bar{p}}^{d}\left(  n\right)  ;
\label{eqLoc.2}%
\end{equation}
if $p=p_{1}+p_{2}$ with $\deg p_{i}\leq\deg p$ ($i=1,2$), then%
\begin{equation}
M_{p}^{d}\left(  n\right)  =\left[  M_{p_{1}}^{d}\left(  n\right)  \right]
_{\eta}+\left[  M_{p_{2}}^{d}\left(  n\right)  \right]  _{\eta}.
\label{eqLoc.3}%
\end{equation}

The main result of this section (Theorem \ref{ThmLoc.2} below) provides a
concrete description of $M_{p}^{d}\left(  n\right)  $ as a linear combination
of certain compressions of $M^{d}\left(  n\right)  $ corresponding to the
monomial terms of $p$. In order to state this result, we require a preliminary
lemma and some additional notation.

\begin{lemma}
\label{LemLoc.1}For $r,s\in\mathbb{Z}_{+}^{d}$ with $\left\vert r\right\vert
+\left\vert s\right\vert \leq2k$, there exist $i,j\in\mathbb{Z}_{+}^{d}$ such
that%
\[
\bar{z}^{r}z^{s}=\bar{z}^{i}z^{j}\bar{z}^{r-i}z^{s-j}\text{\qquad and\qquad
}\left\vert i\right\vert +\left\vert j\right\vert ,\left\vert r-i\right\vert
+\left\vert s-j\right\vert \leq k.
\]
\end{lemma}

\begin{proof}
\textit{Case (i}): $\left|  r\right|  ,\left|  s\right|  \leq k$; let $i=r$,
$j=0$. \textit{Case (ii}): $k<\left|  r\right|  $. We have $r=\left(
r_{1},\dots,r_{d}\right)  $ with $\left|  r\right|  =r_{1}+\dots+r_{d}>k$.
Choose $r^{\prime}\equiv\left(  r_{i}^{\prime},\dots,r_{d}^{\prime}\right)
\in\mathbb{Z}_{+}^{d}$ so that $0\leq r_{i}^{\prime}\leq r_{i}$ ($1\leq i\leq
d$) and $\left|  r^{\prime}\right|  =r_{1}^{\prime}+\dots+r_{d}^{\prime}=k$.
With $i=r^{\prime}$, $j=0$, we have $\left|  r-i\right|  +\left|  s-j\right|
=\left|  r-r^{\prime}\right|  +\left|  s\right|  =\left(  r_{1}-r_{1}^{\prime
}\right)  +\dots+\left(  r_{d}-r_{d}^{\prime}\right)  +\left|  s\right|
=\left|  r\right|  +\left|  s\right|  -\left|  r^{\prime}\right|  \leq2k-k=k$.
\textit{Case (iii}): $k<\left|  s\right|  $; similar to Case (ii).
\end{proof}

For $p\left(  z,\bar{z}\right)  $ as above (with $\delta\equiv\deg p=2k$ or
$2k-1$), we write $p(z,\bar{z})\equiv\sum\limits_{r,s\in\mathbb{Z}_{+}%
^{d},\;\left|  r\right|  +\left|  s\right|  \leq\delta}a_{rs}\bar{z}^{r}z^{s}%
$. Lemma \ref{LemLoc.1} shows that for each $r,s\in\mathbb{Z}_{+}^{d}$ with
$\left|  r\right|  +\left|  s\right|  \leq\delta$, there are tuples $i\equiv
i\left(  r,s,k\right)  $, $j\equiv j\left(  r,s,k\right)  $, $t\equiv t\left(
r,s,k\right)  $, $u\equiv u\left(  r,s,k\right)  $ in $\mathbb{Z}_{+}^{d}$,
such that $i+t=r$, $j+u=s$, $\left|  i\right|  +\left|  j\right|  ,\left|
t\right|  +\left|  u\right|  \leq k$. In the sequel, $_{\left[  \bar{Z}%
^{u}Z^{t};1,\eta\right]  }M^{d}\left(  n\right)  _{\left[  \bar{Z}^{i}%
Z^{j};1,\eta\right]  }$ denotes the compression of $M^{d}\left(  n\right)  $
to the first $\eta$ rows that are indexed by multiples of $\bar{Z}^{u}Z^{t}$
and to the first $\eta$ columns that are indexed by multiples of $\bar{Z}%
^{i}Z^{j}$.

\begin{theorem}
\label{ThmLoc.2}$M_{p}^{d}\left(  n\right)  =\sum\limits_{\substack{r,s\in
\mathbb{Z}_{+}^{d}\\\left|  r\right|  +\left|  s\right|  \leq\delta
}}a_{rs\left[  \bar{Z}^{u}Z^{t};1,\eta\right]  }M^{d}\left(  n\right)
_{\left[  \bar{Z}^{i}Z^{j};1,\eta\right]  }$.
\end{theorem}

For the proof of Theorem \ref{ThmLoc.2}, we require several preliminary
results. Let $0\leq k\leq n$ and let $r,s\in\mathbb{Z}_{+}^{d}$, with $\left|
r\right|  +\left|  s\right|  \leq2k$. From Lemma \ref{LemLoc.1}, we have
$i,j\in\mathbb{Z}_{+}^{d}$ with $\left|  i\right|  +\left|  j\right|  ,\left|
r-i\right|  +\left|  s-j\right|  \leq k$.

\begin{lemma}
\label{LemLoc.3}For $f,g\in\mathbb{C}_{n-k}^{d}\left[  z,\bar{z}\right]  $,%
\[
\left\langle M_{\bar{z}^{r}z^{s}}^{d}\left(  n\right)  \hat{f}%
,\hat{g}\right\rangle =\left\langle M^{d}\left(  n\right)  \left(  \bar{z}%
^{i}z^{j}f\right)  \symbol{94},\left(  z^{r-i}\bar{z}^{s-j}g\right)
\symbol{94}\,\right\rangle .
\]
\end{lemma}

\begin{proof}
Let $\left|  r\right|  +\left|  s\right|  =2l$ or $2l-1$; if $l<k$,
$M_{\bar{z}^{r}z^{s}}^{d}\left(  n\right)  $ has size $\eta\left(
d,n-l\right)  $; in this case we regard $\mathbb{C}_{n-k}^{d}\left[  z,\bar
{z}\right]  $ as embedded in $\mathbb{C}_{n-l}^{d}\left[  z,\bar{z}\right]  $
and take coefficient vectors $\hat{f}$, $\hat{g}$ relative to $\mathbb{C}%
_{n-l}^{d}\left[  z,\bar{z}\right]  $; in any case, $\bar{z}^{i}z^{j}f$ and
$z^{r-i}\bar{z}^{s-j}g$ are elements of $\mathbb{C}_{n}\left[  z,\bar
{z}\right]  $, so $\left(  \bar{z}^{i}z^{j}f\right)  \symbol{94}$ and $\left(
z^{r-i}\bar{z}^{s-j}g\right)  \symbol{94}$ are computed relative to
$\mathbb{C}_{n}\left[  z,\bar{z}\right]  $. We have%
\begin{align*}
\left\langle M_{\bar{z}^{r}z^{s}}^{d}\left(  n\right)  \hat{f}%
,\hat{g}\right\rangle  &  =\Lambda\left(  \bar{z}^{r}z^{s}f\bar{g}\right) \\
&  =\Lambda\left(  \bar{z}^{i}z^{j}f\cdot\left(  z^{r-i}\bar{z}^{s-j}g\right)
^{-}\right) \\
&  =\left\langle M^{d}\left(  n\right)  \left(  \bar{z}^{i}z^{j}f\right)
\symbol{94},\left(  \bar{z}^{s-j}z^{r-i}g\right)  \symbol{94}\,\right\rangle .
\end{align*}
\end{proof}

\begin{proposition}
\label{ProLoc.4}Let $0\leq k\leq n$. Let $r,s,t,u,q,v\in\mathbb{Z}_{+}^{d}$
satisfy $\left\vert r\right\vert +\left\vert s\right\vert \leq2k$, $\left\vert
t\right\vert +\left\vert u\right\vert ,\left\vert q\right\vert +\left\vert
v\right\vert \leq n-k$. Then%
\[
\left\langle M_{\bar{z}^{r}z^{s}}^{d}\left(  n\right)  \widehat{\bar{z}%
^{q}z^{v}},\widehat{\bar{z}^{t}z^{u}}\right\rangle =\gamma_{r+q+u,s+v+t}.
\]
\end{proposition}

\begin{proof}
From Lemma \ref{LemLoc.1}, we have $i,j\in\mathbb{Z}_{+}^{d}$ such that
$\left|  i\right|  +\left|  j\right|  ,\left|  r-i\right|  +\left|
s-j\right|  \leq k$. Lemma \ref{LemLoc.3} implies that
\begin{align*}
\left\langle M_{\bar{z}^{r}z^{s}}^{d}\left(  n\right)  \widehat{\bar{z}%
^{q}z^{v}},\widehat{\bar{z}^{t}z^{u}}\right\rangle  &  =\left\langle
M^{d}\left(  n\right)  \left(  \bar{z}^{i}z^{j}\bar{z}^{q}z^{v}\right)
\symbol{94},\left(  z^{r-i}\bar{z}^{s-j}\bar{z}^{t}z^{u}\right)
\symbol{94}\,\right\rangle \\
&  =\left\langle M^{d}\left(  n\right)  \left(  \bar{z}^{i+q}z^{j+v}\right)
\symbol{94},\left(  \bar{z}^{s+t-j}z^{r+u-i}\right)  \symbol{94}%
\,\right\rangle \\
&  =\gamma_{\left(  i+q\right)  +\left(  r+u-i\right)  ,\left(  j+v\right)
+\left(  s+t-j\right)  }=\gamma_{q+r+u,s+v+t}.
\end{align*}
\end{proof}

\begin{lemma}
\label{LemLoc.5}Let $0\leq k\leq n$ and let $\eta=\eta\left(  d,n-k\right)  $.
Suppose $p,q,l,m\in\mathbb{Z}_{+}^{d}$ satisfy $\left|  p\right|  +\left|
q\right|  ,\left|  l\right|  +\left|  m\right|  \leq k$ and set
\[
M:=_{\left[  \bar{Z}^{m}Z^{l};1,\eta\right]  }M^{d}\left(  n\right)  _{\left[
\bar{Z}^{p}Z^{q};1,\eta\right]  }.
\]
Then $M=\left[  M_{\bar{z}^{p}z^{q}\cdot\bar{z}^{l}z^{m}}^{d}\left(  n\right)
\right]  _{\eta}$, the compression of $M_{\bar{z}^{p}z^{q}\cdot\bar{z}%
^{l}z^{m}}^{d}$ to its first $\eta$ rows and columns.
\end{lemma}

\begin{proof}
The columns of $M$ are indexed by $\bar{Z}^{p+i}Z^{q+j}$, $i,j\in
\mathbb{Z}_{+}^{d}$, $\left|  i\right|  +\left|  j\right|  \leq n-k$, and the
rows are indexed by $\bar{Z}^{m+a}Z^{l+b}$, $a,b\in\mathbb{Z}_{+}^{d}$,
$\left|  a\right|  +\left|  b\right|  \leq n-k$. The entry in row $\bar
{Z}^{m+a}Z^{l+b}$, column $\bar{Z}^{p+i}Z^{q+j}$ of $M$ is thus
\[
\left\langle M^{d}\left(  n\right)  \left(  \bar{z}^{p+i}z^{q+j}\right)
\symbol{94},\left(  \bar{z}^{m+a}z^{l+b}\right)  \symbol{94}\,\right\rangle
=\gamma_{p+i+l+b,q+j+m+a}.
\]
The corresponding entry of $M_{\bar{z}^{p}z^{q}\cdot\bar{z}^{l}z^{m}}%
^{d}\left(  n\right)  $, in row $\bar{Z}^{a}Z^{b}$, column $\bar{Z}^{i}Z^{j}$,
is $\left\langle M_{\bar{z}^{p}z^{q}\cdot\bar{z}^{l}z^{m}}^{d}\left(
n\right)  \widehat{\bar{z}^{i}z^{j}},\widehat{\bar{z}^{a}z^{b}}\right\rangle
$, which, by Proposition \ref{ProLoc.4}, is also equal to $\gamma
_{p+i+l+b,q+j+m+a}$.
\end{proof}

\begin{proof}
[Proof of Theorem \textup{\ref{ThmLoc.2}}]We have $1\leq k\leq n$ and
\[
p\equiv p\left(  z,\bar{z}\right)  =\sum\limits_{r,s\in\mathbb{Z}_{+}%
^{d},\;\left|  r\right|  +\left|  s\right|  \leq\delta}a_{rs}\bar{z}^{r}%
z^{s},
\]
with $\delta\equiv\deg p$ ($=2k$ or $2k-1$). The size of $M_{p}^{d}\left(
n\right)  $ is thus $\eta\times\eta$, where $\eta=\eta\left(  d,n-k\right)  $.
By (\ref{eqLoc.3}) and the uniqueness of $M_{p}^{d}\left(  n\right)  $, we
have
\begin{equation}
M_{p}^{d}\left(  n\right)  =\sum_{r,s\in\mathbb{Z}_{+}^{d},\;\left|  r\right|
+\left|  s\right|  \leq\delta}a_{rs}\left[  M_{\bar{z}^{r}z^{s}}^{d}\left(
n\right)  \right]  _{\eta}. \label{eqLoc.2bis}%
\end{equation}
From Lemma \ref{LemLoc.1}, for each $r,s\in\mathbb{Z}_{+}^{d}$ with $\left|
r\right|  +\left|  s\right|  \leq\delta$, we have $i\equiv i\left(
r,s,k\right)  ,j\equiv j\left(  r,s,k\right)  ,t\equiv t\left(  r,s,k\right)
,u\equiv u\left(  r,s,k\right)  \in\mathbb{Z}_{+}^{d}$ with $i+t=r$, $j+u=s$,
$\left|  i\right|  +\left|  j\right|  ,\left|  t\right|  +\left|  u\right|
\leq k$. Lemma \ref{LemLoc.5} implies that for each $r$, $s$,%
\begin{align*}
\left[  M_{\bar{z}^{r}z^{s}}^{d}\left(  n\right)  \right]  _{\eta}  &
=\left[  M_{\bar{z}^{i}z^{j}\cdot\bar{z}^{t}z^{u}}^{d}\left(  n\right)
\right]  _{\eta}\\
&  =_{\left[  \bar{Z}^{u}Z^{t};1,\eta\right]  }M^{d}\left(  n\right)
_{\left[  \bar{Z}^{i}Z^{j};1,\eta\right]  },
\end{align*}
so the result follows from (\ref{eqLoc.2bis}).
\end{proof}

We conclude this section with an analogue of Theorem \ref{ThmLoc.2} for real
moment matrices. \ Given a real moment matrix $\mathcal{M}^{N}(n)\equiv
\mathcal{M}^{N}(n)(\beta)$, let $k\leq n$, and let $p\in\mathbb{C}%
[t_{1},...,t_{N}]$, with $\deg p=2k$ or $2k-1$. \ The \textit{localizing
matrix} $\mathcal{M}_{p}^{N}(n)$ has size $\tau\equiv\tau(N,n-k):=\left(
\begin{array}
[c]{c}%
n-k+N\\
N
\end{array}
\right)  $ and is uniquely determined by
\begin{equation}
\left\langle \mathcal{M}_{p}^{N}(n)\hat{f},\hat{g}\right\rangle =\Lambda
_{\beta}(pf\bar{g})\;\;(f,g\in\mathbb{C}_{n-k}^{N}[t]); \label{neweq}%
\end{equation}
\ if $\beta$ has a representing measure supported in $K_{p}:=\{t\in
\mathbb{R}^{N}:p(t)\geq0\}$, then clearly $\mathcal{M}_{p}^{N}(n)\geq0$.
\ Write $p(t)\equiv\sum_{i\in\mathbb{Z}_{+}^{n},\left|  i\right|  \leq\deg
p}a_{i}t^{i}$. \ For each $i$, there exist (non-unique) $r\equiv r(i)$,
$s\equiv s(i)$ in $\mathbb{Z}_{+}^{n}$ such that $r+s=i$ and $\left|
r\right|  ,\left|  s\right|  \leq k$; thus $p(t)=\sum_{i\in\mathbb{Z}_{+}%
^{n},\left|  i\right|  \leq\deg p}a_{i}t^{r(i)}t^{s(i)}$. \ Let $_{\left[
T^{r};1,\tau\right]  }\mathcal{M}^{N}(n)_{\left[  T^{s};1,\tau\right]  }$
denote the compression of $\mathcal{M}^{N}(n)$ to the first $\tau$ rows that
are indexed by multiples of $T^{r}$ and to the first $\tau$ columns that are
indexed by multiples of $T^{s}$.

\begin{theorem}
\label{thm3.6}$\mathcal{M}_{p}^{N}(n)=\sum_{i\in\mathbb{Z}_{+}^{n},\left|
i\right|  \leq\deg p}a_{i\;\left[  T^{r};1,\tau\right]  }\mathcal{M}%
^{N}(n)_{\left[  T^{s};1,\tau\right]  }$.
\end{theorem}

The proof of Theorem \ref{thm3.6} follows by formal repetition of the proof of
Theorem \ref{ThmLoc.2}; we omit the details. \ Example \ref{IntroExample1}
illustrates Theorem \ref{thm3.6} with $N=3,n=1,\deg p=2$.

\section{\label{Fla}Flat extensions of positive localizing matrices}

In this section we present a flat extension theorem for positive localizing
matrices, which provides the main tool for proving Theorem \ref{IntroTheorem1}%
. Suppose $M^{d}\left(  n\right)  \left(  \gamma\right)  $ is positive and
admits a flat extension $M^{d}\left(  n+1\right)  $; thus, there is a matrix
$W$ such that $M^{d}\left(  n+1\right)  $ admits a block decomposition of the
form%
\begin{equation}
M^{d}\left(  n+1\right)  =\left(
\begin{array}
[c]{cc}%
M^{d}\left(  n\right)  & B^{d}\left(  n+1\right) \\
B^{d}\left(  n+1\right)  ^{\ast} & C^{d}\left(  n+1\right)
\end{array}
\right)  , \label{eqFla.1}%
\end{equation}
where $B^{d}\left(  n+1\right)  =M^{d}\left(  n\right)  W$ and $C^{d}\left(
n+1\right)  =W^{\ast}M^{d}\left(  n\right)  W$. It follows from Theorem
\ref{MomThm2} that $M^{d}\left(  n+1\right)  $ admits a unique positive flat
extension $M^{d}(\infty)$ and that $M^{d}(\infty)$ admits a representing
measure $\nu$. \ In particular, $M^{d}(n+1)\equiv M^{d}(n+1)[\nu]$ is positive
and recursively generated, and $M^{d}\left(  n+1\right)  $ admits unique
successive positive, flat moment matrix extensions $M^{d}\left(  n+2\right)
\equiv M^{d}(n+2)[\nu],M^{d}\left(  n+3\right)  \equiv M^{d}(n+3)[\nu],\dots$.
Thus, if $p\in\mathbb{C}^{d}\left[  z,\bar{z}\right]  $ and $k:=\left[
\left(  1+\deg p\right)  /2\right]  \leq n$, we may consider $M_{p}^{d}\left(
n+k\right)  $ and $M_{p}^{d}\left(  n+k+1\right)  $.

\begin{theorem}
\label{ThmFla.1}Suppose $M^{d}\left(  n\right)  \left(  \gamma\right)  \geq0$
admits a flat extension%
\[
M^{d}\left(  n+1\right)  =\left(
\begin{array}
[c]{cc}%
M^{d}\left(  n\right)  & M^{d}\left(  n\right)  W\\
W^{\ast}M^{d}\left(  n\right)  & W^{\ast}M^{d}\left(  n\right)  W
\end{array}
\right)  .
\]
Let $p\in\mathbb{C}^{d}\left[  z,\bar{z}\right]  $, with $\deg p=2k$ or
$2k-1$. If $M_{p}^{d}\left(  n+k\right)  \geq0$, then%
\begin{equation}
M_{p}^{d}\left(  n+k+1\right)  =\left(
\begin{array}
[c]{cc}%
M_{p}^{d}\left(  n+k\right)  & M_{p}^{d}\left(  n+k\right)  W\\
W^{\ast}M_{p}^{d}\left(  n+k\right)  & W^{\ast}M_{p}^{d}\left(  n+k\right)  W
\end{array}
\right)  ; \label{eqFla.2}%
\end{equation}
in particular, $M_{p}^{d}\left(  n+k+1\right)  $ is a flat, positive extension
of $M_{p}\left(  n+k\right)  $.
\end{theorem}

\begin{remark}
\label{RemFla.0}In Theorem \textup{\ref{ThmFla.1}}, we are not assuming that
$M_{p}^{d}\left(  n+k\right)  $ is a moment matrix; rather, in Section
\textup{\ref{Exi}} we will prove that under the hypotheses of Theorem
\textup{\ref{ThmFla.1}}, both $M_{p}^{d}\left(  n+k\right)  $ and $M_{p}%
^{d}\left(  n+k+1\right)  $ are indeed moment matrices.
\end{remark}

The proof of Theorem \ref{ThmFla.1} is based on a computational description of
$M_{p}^{d}\left(  n+k+1\right)  $, and to derive this we require some
additional notation. For $m>0$, let $A$ be a matrix of size $\eta\left(
d,m\right)  $ with rows and columns $\left\{  \bar{Z}^{a}Z^{b}\right\}
_{a,b\in\mathbb{Z}_{+}^{d},\;\left|  a\right|  +\left|  b\right|  \leq m}$,
ordered lexicographically. Suppose $i,j\in\mathbb{Z}_{+}^{d}$, with $\left|
i\right|  +\left|  j\right|  \leq m$, and suppose there are at least $\beta$
columns of $A$ that are indexed by multiples of $\bar{Z}^{i}Z^{j}$. Suppose
$u,t\in\mathbb{Z}_{+}^{d}$, with $\left|  u\right|  +\left|  t\right|  \leq
m$, and suppose there are at least $\upsilon$ rows of $A$ that are indexed by
multiples of $\bar{Z}^{u}Z^{t}$. For $\alpha\leq\beta$ and $\rho\leq\upsilon$,
let $_{\left[  \bar{Z}^{u}Z^{t},\rho,\upsilon\right]  }A_{\left[  \bar{Z}%
^{i}Z^{j},\alpha,\beta\right]  }$ denote the compression of $A$ to the
$\alpha$-th through $\beta$-th consecutive columns indexed by multiples of
$\bar{Z}^{i}Z^{j}$ and to the $\rho$-th through $\upsilon$-th consecutive rows
indexed by multiples of $\bar{Z}^{u}Z^{t}$. We omit the proof of the following
elementary result.

\begin{lemma}
\label{LemFla.2}$\left(  _{\left[  \bar{Z}^{u}Z^{t};\rho,\upsilon\right]
}A_{\left[  \bar{Z}^{i}Z^{j};\alpha,\beta\right]  }\right)  ^{\ast}=_{\left[
\bar{Z}^{i}Z^{j};\alpha,\beta\right]  }(A^{\ast})_{\left[  \bar{Z}^{u}%
Z^{t};\rho,\upsilon\right]  }$.
\end{lemma}

(Here, the convention is that rows and columns of $A^{\ast}$ are indexed in
the same way as the rows and columns of $A$, as $\left\{  \bar{Z}^{a}%
Z^{b}\right\}  _{a,b\in\mathbb{Z}_{+}^{d},\;\left\vert a\right\vert
+\left\vert b\right\vert \leq m}$.)

To prove Theorem \ref{ThmFla.1}, we will first obtain analogues of
(\ref{eqFla.2}) for each monomial term of $p$. To this end, let $\delta=\deg
p$ ($=2k$ or $2k-1$); write $p$ as $p\left(  z,\bar{z}\right)  \equiv
\sum\limits_{r,s\in\mathbb{Z}_{+}^{d},\;\left|  r\right|  +\left|  s\right|
\leq\delta}a_{rs}\bar{z}^{r}z^{s}$. Recall from Section \ref{Loc} that
\[
\eta_{2}\equiv\operatorname*{size}M_{p}^{d}\left(  n+k+1\right)  =\eta\left(
d,\left(  n+k+1\right)  -k\right)  =\eta\left(  d,n+1\right)
=\operatorname*{size}M^{d}\left(  n+1\right)
\]
and
\[
\eta_{1}\equiv\operatorname*{size}M_{p}^{d}\left(  n+k\right)  =\eta\left(
d,\left(  n+k\right)  -k\right)  =\eta\left(  d,n\right)
=\operatorname*{size}M^{d}\left(  n\right)  .
\]
Let $p_{rs}=\bar{z}^{r}z^{s}$; from Lemma \ref{LemLoc.1}, we have $\bar{z}%
^{r}z^{s}=\bar{z}^{i\left(  r,s,k\right)  }z^{j\left(  r,s,k\right)  }%
\cdot\bar{z}^{t\left(  r,s,k\right)  }z^{u\left(  r,s,k\right)  }$, where
$i\equiv i\left(  r,s,k\right)  ,j\equiv j\left(  r,s,k\right)  ,t\equiv
t\left(  r,s,k\right)  ,$ and $u\equiv u\left(  r,s,k\right)  \in
\mathbb{Z}_{+}^{d}$ satisfy $r=i+t$, $s=j+u$, $\left|  i\right|  +\left|
j\right|  ,\left|  t\right|  +\left|  u\right|  \leq k$. Lemma \ref{LemLoc.5}
(applied with $n$ replaced by $n+k+1$) shows that
\begin{equation}
\left[  M_{p_{rs}}^{d}\left(  n+k+1\right)  \right]  _{\eta_{2}}=_{\left[
\bar{Z}^{u}Z^{t};1,\eta_{2}\right]  }M^{d}\left(  n+k+1\right)  _{\left[
\bar{Z}^{i}Z^{j};1,\eta_{2}\right]  }; \label{eqFla.3}%
\end{equation}
similarly,%
\begin{equation}
\left[  M_{p_{rs}}^{d}\left(  n+k\right)  \right]  _{\eta_{1}}=_{\left[
\bar{Z}^{u}Z^{t};1,\eta_{1}\right]  }M^{d}\left(  n+k\right)  _{\left[
\bar{Z}^{i}Z^{j};1,\eta_{1}\right]  }. \label{eqFla.4}%
\end{equation}

We next use (\ref{eqFla.3}) and (\ref{eqFla.4}) to relate $\left[  M_{p_{rs}%
}^{d}\left(  n+k+1\right)  \right]  _{\eta_{2}}$ to $\left[  M_{p_{rs}}%
^{d}\left(  n+k\right)  \right]  _{\eta_{1}}$ via a block decomposition of
$\left[  M_{p_{rs}}^{d}\left(  n+k+1\right)  \right]  _{\eta_{2}}$. From
(\ref{eqFla.3}), note that the columns of $\left[  M_{p_{rs}}^{d}\left(
n+k+1\right)  \right]  _{\eta_{2}}$ are compressions of the first $\eta_{2}$
columns of $M^{d}\left(  n+k+1\right)  $ that are indexed by multiples of
$\bar{Z}^{i}Z^{j}$; these monomials are ordered as $\left\{  \bar{Z}^{i+i_{q}%
}Z^{j+j_{q}}\right\}  _{q=1}^{\eta_{2}}$, where $\left\{  \bar{Z}^{i_{q}%
}Z^{j_{q}}\right\}  _{q=1}^{\eta_{2}}$ is the lexicographic ordering of the
first $\eta_{2}$ monomials in $\mathbb{C}^{d}\left[  z,\bar{z}\right]  $. In
particular, from (\ref{eqFla.4}) we see that the first $\eta_{1}$ of these
monomials also index the columns of $\left[  M_{p_{rs}}^{d}\left(  n+k\right)
\right]  _{\eta_{1}}$. Similarly, the rows of $\left[  M_{p_{rs}}^{d}\left(
n+k+1\right)  \right]  _{\eta_{2}}$ are compressions of rows of \newline
$M^{d}\left(  n+k+1\right)  $ that are indexed by the sequence $\left\{
\bar{Z}^{u+i_{q}}Z^{t+j_{q}}\right\}  _{q=1}^{\eta_{2}}$, and the first
$\eta_{1}$ of these also index the rows of $\left[  M_{p_{rs}}^{d}\left(
n+k\right)  \right]  _{\eta_{1}}$. Further, from (\ref{eqFla.1}) and the above
remarks, it is clear that $\left[  M_{p_{rs}}^{d}\left(  n+k+1\right)
\right]  _{\eta_{1}}=\left[  M_{p_{rs}}^{d}\left(  n+k\right)  \right]
_{\eta_{1}}$. It now follows from the preceding observations that
\newline $\left[  M_{p_{rs}}^{d}\left(  n+k+1\right)  \right]  _{\eta_{2}}$
admits a block decomposition of the form%
\begin{equation}
M_{p_{rs}}^{d}\left(  n+k+1\right)  _{\eta_{2}}=\left(
\begin{array}
[c]{cc}%
\left[  M_{p_{rs}}^{d}\left(  n+k\right)  \right]  _{\eta_{1}} & B_{p_{rs}%
}^{d}\left(  n+k+1\right) \\
D_{p_{rs}}^{d}\left(  n+k+1\right)  & C_{p_{rs}}^{d}\left(  n+k+1\right)
\end{array}
\right)  , \label{eqFla.5}%
\end{equation}
where%
\begin{align}
\left[  M_{p_{rs}}^{d}\left(  n+k\right)  \right]  _{\eta_{1}}  &  =_{\left[
\bar{Z}^{u}Z^{t};1,\eta_{1}\right]  }M^{d}\left(  n+k\right)  _{\left[
\bar{Z}^{i}Z^{j};1,\eta_{1}\right]  },\label{eqFla.6}\\
B_{p_{rs}}^{d}\left(  n+k+1\right)   &  =_{\left[  \bar{Z}^{u}Z^{t};1,\eta
_{1}\right]  }M^{d}\left(  n+k+1\right)  _{\left[  \bar{Z}^{i}Z^{j};1+\eta
_{1},\eta_{2}\right]  },\label{eqFla.7}\\
D_{p_{rs}}^{d}\left(  n+k+1\right)   &  =_{\left[  \bar{Z}^{u}Z^{t};\eta
_{1}+1,\eta_{2}\right]  }M^{d}\left(  n+k+1\right)  _{\left[  \bar{Z}^{i}%
Z^{j};1,\eta_{1}\right]  },\label{eqFla.8}\\
C_{p_{rs}}^{d}\left(  n+k+1\right)   &  =_{\left[  \bar{Z}^{u}Z^{t};\eta
_{1}+1,\eta_{2}\right]  }M^{d}\left(  n+k+1\right)  _{\left[  \bar{Z}^{i}%
Z^{j};\eta_{1}+1,\eta_{2}\right]  }. \label{eqFla.9}%
\end{align}

The following lemma is the first step toward proving an analogue of Theorem
\ref{ThmFla.1} for $p_{rs}$.

\begin{lemma}
\label{LemFla.3}For each $r,s\in\mathbb{Z}_{+}^{d}$ with $\left|  r\right|
+\left|  s\right|  \leq\delta$,%
\[
\left(
\begin{array}
[c]{c}%
\left[  M_{p_{rs}}^{d}\left(  n+k\right)  \right]  _{\eta_{1}}\\
D_{p_{rs}}^{d}\left(  n+k+1\right)
\end{array}
\right)  W=\left(
\begin{array}
[c]{c}%
B_{p_{rs}}^{d}\left(  n+k+1\right) \\
C_{p_{rs}}^{d}\left(  n+k+1\right)
\end{array}
\right)  .
\]
\end{lemma}

\begin{proof}
For $1\leq m\leq\eta_{2}-\eta_{1}$, the $m$-th column of $\left(
\begin{array}
[c]{c}%
B^{d}\left(  n+1\right) \\
C^{d}\left(  n+1\right)
\end{array}
\right)  $ is the $\left(  \eta_{1}+m\right)  $-th column of $M\left(
n+1\right)  $, and is thus of the form $\bar{Z}^{e}Z^{f}\in\mathcal{C}%
_{M\left(  n+1\right)  }$, with $\left|  e\right|  +\left|  f\right|  =n+1$.
If $\left(  \alpha_{a,b}^{\left(  m\right)  }\right)  _{a,b\in\mathbb{Z}%
_{+}^{d},\;\left|  a\right|  +\left|  b\right|  \leq n}$ denotes the $m$-th
column of $W$, then we have%
\begin{equation}
\bar{Z}^{e}Z^{f}=\sum_{\left|  a\right|  +\left|  b\right|  \leq n}%
\alpha_{a,b}^{\left(  m\right)  }\bar{Z}^{a}Z^{b}. \label{eqFla.10}%
\end{equation}
Let $\left\{  V_{a,b}\left(  r,s\right)  \right\}  _{\left|  a\right|
+\left|  b\right|  \leq n}$ denote the lexicographic ordering of the columns
of $\left(
\begin{array}
[c]{c}%
\left[  M_{p_{rs}}^{d}\left(  n+k\right)  \right]  _{\eta_{1}}\\
D_{p_{rs}}^{d}\left(  n+k+1\right)
\end{array}
\right)  $ and let $U_{m}\left(  r,s\right)  $ denote the $m$-th column of
\newline $\left(
\begin{array}
[c]{c}%
B_{p_{rs}}^{d}\left(  n+k+1\right) \\
C_{p_{rs}}^{d}\left(  n+k+1\right)
\end{array}
\right)  $. It suffices to show that
\begin{equation}
U_{m}\left(  r,s\right)  =\sum_{\left|  a\right|  +\left|  b\right|  \leq
n}\alpha_{a,b}^{\left(  m\right)  }V_{a,b}\left(  r,s\right)  .
\label{eqFla.11}%
\end{equation}

Since $M^{d}\left(  n+k+1\right)  $ is a flat, hence positive, extension of
$M^{d}\left(  n+1\right)  $, (\ref{eqFla.10}) also holds in $\mathcal{C}%
_{M^{d}\left(  n+k+1\right)  }$. Now $U_{m}\left(  r,s\right)  $ is the
$\left(  \eta_{1}+m\right)  $-th column of $\left[  M_{p_{rs}}^{d}\left(
n+k+1\right)  \right]  _{\eta_{2}}$, and is thus indexed by the $\left(
\eta_{1}+m\right)  $-th multiple of $\bar{Z}^{i\left(  r,s,k\right)
}Z^{j\left(  r,s,k\right)  }$; thus $U_{m}\left(  r,s\right)  $ is indexed by
$\bar{Z}^{e+i\left(  r,s,k\right)  }Z^{f+j\left(  r,s,k\right)  }$. Since
$M^{d}\left(  n+k+1\right)  $ is recursively generated, (\ref{eqFla.10})
implies that in $\mathcal{C}_{M^{d}\left(  n+k+1\right)  }$ we have%
\begin{equation}
\bar{Z}^{e+i\left(  r,s,k\right)  }Z^{f+j\left(  r,s,k\right)  }=\sum_{\left|
a\right|  +\left|  b\right|  \leq n}\alpha_{a,b}^{\left(  m\right)  }\bar
{Z}^{a+i\left(  r,s,k\right)  }Z^{b+j\left(  r,s,k\right)  }; \label{eqFla.12}%
\end{equation}
thus, via compression of these columns to rows indexed by the first $\eta_{2}$
multiples of $\bar{Z}^{u}Z^{t}$, it follows that the relation in
(\ref{eqFla.12}) holds as well in the column space of $\left[  M_{p_{rs}}%
^{d}\left(  n+k+1\right)  \right]  _{\eta_{2}}$. Since the compression of
$\bar{Z}^{e+i\left(  r,s,k\right)  }Z^{f+j\left(  r,s,k\right)  }$ is
$U_{m}\left(  r,s\right)  $ and the compression of $\bar{Z}^{a+i\left(
r,s,k\right)  }Z^{b+j\left(  r,s,k\right)  }$ is $V_{a,b}\left(  r,s\right)
$, we obtain (\ref{eqFla.11}), so the result follows.
\end{proof}

\begin{lemma}
\label{LemFla.4}For each $r,s\in\mathbb{Z}_{+}^{d}$ with $\left|  r\right|
+\left|  s\right|  \leq\delta$,
\[
D_{p_{rs}}^{d}\left(  n+k+1\right)  =B_{\bar{p}_{rs}}^{d}\left(  n+k+1\right)
^{\ast}=W^{\ast}\left[  M_{p_{rs}}^{d}\left(  n+k\right)  \right]  _{\eta_{1}%
}.
\]
\end{lemma}

\begin{proof}
Applying Lemma \ref{LemFla.3} to $\bar{p}_{rs}$ ($\equiv\bar{Z}^{s}Z^{r}$), we
have%
\begin{align*}
B_{\bar{p}_{rs}}^{d}\left(  n+k+1\right)   &  =\left[  M_{\bar{p}_{rs}}%
^{d}\left(  n+k\right)  \right]  _{\eta_{1}}W\\
&  =\left[  M_{p_{rs}}^{d}\left(  n+k\right)  ^{\ast}\right]  _{\eta_{1}%
}W\;\;\;\text{(by (\ref{eqLoc.2}))}\\
&  =\left[  M_{p_{rs}}^{d}\left(  n+k\right)  \right]  _{\eta_{1}}^{\ast}W,
\end{align*}
whence $B_{\bar{p}_{rs}}^{d}\left(  n+k+1\right)  ^{\ast}=W^{\ast}\left[
M_{p_{rs}}^{d}\left(  n+k\right)  \right]  _{\eta_{1}}$. Now
\begin{align*}
D_{p_{rs}}^{d}\left(  n+k+1\right)  ^{\ast}  &  =_{\left[  \bar{Z}^{i}%
Z^{j};1,\eta_{1}\right]  }M^{d}\left(  n+k+1\right)  _{\left[  \bar{Z}%
^{u}Z^{t};\eta_{1}+1,\eta_{2}\right]  }^{\ast}\\
&  \qquad\qquad\qquad\qquad\qquad\qquad\text{(Lemma \ref{LemFla.2})}\\
&  =_{\left[  \bar{Z}^{i}Z^{j};1,\eta_{1}\right]  }M^{d}\left(  n+k+1\right)
_{\left[  \bar{Z}^{u}Z^{t};\eta_{1}+1,\eta_{2}\right]  }\\
&  \qquad\qquad\qquad\qquad\qquad\qquad\text{(since }M^{d}(n+k+1)\geq
0\text{)}\\
&  =B_{\bar{Z}^{u}Z^{t}\cdot Z^{i}\bar{Z}^{j}}^{d}\left(  n+k+1\right)
\;\;\;\text{(by (\ref{eqFla.7}) applied to }\bar{Z}^{s}Z^{r}\text{)}\\
&  =B_{\bar{p}_{rs}}^{d}\left(  n+k+1\right)  ;
\end{align*}
thus $D_{p_{rs}}^{d}\left(  n+k+1\right)  =B_{\bar{p}_{rs}}^{d}\left(
n+k+1\right)  ^{\ast}$ and the proof is complete.
\end{proof}

\begin{proof}
[Proof of Theorem \textup{\ref{ThmFla.1}}]Using the uniqueness of $M_{p}%
^{d}\left(  n+k\right)  $ and of \newline $M_{p}^{d}\left(  n+k+1\right)  $,
it follows from (\ref{eqFla.5})--(\ref{eqFla.9}) that $M_{p}^{d}\left(
n+k+1\right)  $ admits a block decomposition of the form%
\[
M_{p}^{d}\left(  n+k+1\right)  =\left(
\begin{array}
[c]{cc}%
M_{p}^{d}\left(  n+k\right)  & B_{p}^{d}\left(  n+k+1\right) \\
D_{p}^{d}\left(  n+k+1\right)  & C_{p}^{d}\left(  n+k+1\right)
\end{array}
\right)  ,
\]
where
\begin{align*}
M_{p}^{d}\left(  n+k\right)   &  =\sum_{\left|  r\right|  +\left|  s\right|
\leq\delta}\alpha_{rs}\left[  M_{p_{rs}}^{d}\left(  n+k\right)  \right]
_{\eta_{1}}\,,\\
B_{p}^{d}\left(  n+k+1\right)   &  =\sum_{\left|  r\right|  +\left|  s\right|
\leq\delta}\alpha_{rs}B_{p_{rs}}^{d}\left(  n+k+1\right)  ,\\
D_{p}^{d}\left(  n+k+1\right)   &  =\sum_{\left|  r\right|  +\left|  s\right|
\leq\delta}\alpha_{rs}D_{p_{rs}}^{d}\left(  n+k+1\right)  ,\\
C_{p}^{d}\left(  n+k+1\right)   &  =\sum_{\left|  r\right|  +\left|  s\right|
\leq\delta}\alpha_{rs}C_{p_{rs}}^{d}\left(  n+k+1\right)  .
\end{align*}
Lemma \ref{LemFla.3} implies $B_{p}^{d}\left(  n+k+1\right)  =M_{p}^{d}\left(
n+k\right)  W$, and Lemma \ref{LemFla.4} implies $D_{p}^{d}\left(
n+k+1\right)  =W^{\ast}M_{p}^{d}\left(  n+k\right)  $. Now Lemmas
\ref{LemFla.3} and \ref{LemFla.4} imply $C_{p}^{d}\left(  n+k+1\right)
=D_{p}^{d}\left(  n+k+1\right)  W=W^{\ast}M_{p}^{d}\left(  n+k\right)  W$,
whence (\ref{eqFla.2}) holds. Since $M_{p}^{d}\left(  n+k\right)  $ is
positive, (\ref{eqFla.2}) implies that $M_{p}^{d}\left(  n+k+1\right)  $ is
positive and that $\operatorname{rank}M_{p}^{d}\left(  n+k+1\right)
=\operatorname{rank}M_{p}^{d}\left(  n+k\right)  $ \ (cf. \cite{tcmp4}).
\end{proof}

\section{\label{Exi}Existence of minimal representing measures supported in
semi-algebraic sets}

We begin with the analogue of Theorem \ref{IntroTheorem1} for the truncated
complex multivariable $K$-moment problem. \ Recall that if $M^{d}\left(
n\right)  (\gamma)\;(\geq0)$ has a flat extension $M^{d}\left(  n+1\right)  $,
then $M^{d}(n+1)$ admits unique recursive flat (positive) extensions
$M^{d}\left(  n+2\right)  ,M^{d}\left(  n+3\right)  ,...\;$(Theorem
\ref{MomThm2})$\mathbf{.}$

\begin{theorem}
\label{ThmExi.1}Let $\gamma\equiv\gamma^{(2n)}=\{\gamma_{i}\}_{i\in
\mathbb{Z}_{+}^{d},\;\left|  i\right|  \leq2n}$ be a complex sequence and let
$\mathcal{P}\equiv\{p_{i}\}_{i=1}^{m}\subseteq\mathbb{C}^{d}\left[  z,\bar
{z}\right]  $ with $\deg p_{i}=2k_{i}$ or $2k_{i}-1\;(1\leq i\leq m)$. \ Let
$M\equiv M^{d}\left(  n\right)  (\gamma)$ and let $r:=\operatorname{rank}M$.
\ There exists a (minimal) $r$-atomic representing measure for $\gamma$
supported in $K_{\mathcal{P}}$ if and only if $M\geq0$ and $M$ admits a flat
extension $M^{d}\left(  n+1\right)  $ for which $M_{p_{i}}^{d}\left(
n+k_{i}\right)  \geq0\;(1\leq i\leq m)$. \ In this case, $M^{d}(n+1)$ admits a
unique representing measure $\nu$, which is an $r$-atomic (minimal)
$K_{\mathcal{P}}$-representing measure for $\gamma$; moreover, $\nu$ has
precisely $r-\operatorname{rank}M_{p_{i}}^{d}\left(  n+k_{i}\right)  $ atoms
in $\mathcal{Z}\left(  p_{i}\right)  \;(1\leq i\leq m)$.
\end{theorem}

\begin{proof}
Suppose $M^{d}(n)(\gamma)$ is positive and admits a flat extension
$M^{d}\left(  n+1\right)  $ for which $M_{p_{i}}^{d}\left(  n+k_{i}\right)
\geq0\;(1\leq i\leq m)$. \ \cite[Corollary 7.9]{tcmp1} and \cite[Theorem
7.7]{tcmp1} imply that $M^{d}\left(  n+1\right)  $ admits a unique flat
(positive) extension $M^{d}\left(  \infty\right)  $, and that $M^{d}\left(
\infty\right)  $ admits an $r$-atomic representing measure $\nu\equiv
\sum_{j=1}^{r}\rho_{j}\delta_{\omega_{j}}$, with $\rho_{j}>0$ and $\omega
_{j}\in\mathbb{C}^{d}\;(1\leq j\leq r)$. \ Theorem \ref{thm12} implies that
$\nu$ is the unique representing measure for $M^{d}(n+1)$. \ We will show that
$\operatorname*{supp}\nu\subseteq K_{\mathcal{P}}$. \ Fix $i$, $1\leq i\leq
m$. \ Since $M_{p_{i}}^{d}\left(  n+k_{i}\right)  \geq0$, repeated application
of Theorem \ref{ThmFla.1} shows that $M_{p_{i}}^{d}\left(  \infty\right)  $ is
a flat, positive extension of $M_{p_{i}}^{d}\left(  n+k_{i}\right)  $;
moreover,
\begin{equation}
\left\langle M_{p_{i}}^{d}\left(  \infty\right)  \hat{f},\hat{g}\right\rangle
=\int p_{i}\,f\,\bar{g}\,d\nu,\qquad f,g\in\mathbb{C}^{d}\left[  z,\bar
{z}\right]  . \label{eqExi.2}%
\end{equation}
Fix $j$, $1\leq j\leq r$, and let%
\[
f_{j}\left(  z,\bar{z}\right)  =\frac{\left\|  z-\omega_{1}\right\|
^{2}\cdots\left\|  z-\omega_{j-1}\right\|  ^{2}\left\|  z-\omega
_{j+1}\right\|  ^{2}\cdots\left\|  z-\omega_{r}\right\|  ^{2}}{\left\|
\omega_{j}-\omega_{1}\right\|  ^{2}\cdots\left\|  \omega_{j}-\omega
_{j-1}\right\|  ^{2}\left\|  \omega_{j}-\omega_{j+1}\right\|  ^{2}%
\cdots\left\|  \omega_{j}-\omega_{r}\right\|  ^{2}}%
\]
(where, for $z\equiv\left(  z_{1},\dots,z_{d}\right)  $, $\left\|  z\right\|
^{2}:=\sum\bar{z}_{i}z_{i}\in\mathbb{C}^{d}\left[  z,\bar{z}\right]  $)$.$ Now
$f_{j}\in\mathbb{C}^{d}\left[  z,\bar{z}\right]  $, so by (\ref{eqExi.2}),%
\begin{align*}
0\leq\left\langle M_{p_{i}}^{d}\left(  \infty\right)  \hat{f}_{j},\hat{f}%
_{j}\right)   &  =\int p_{i}\left|  f_{j}\right|  ^{2}\,d\nu\\
&  =\sum_{k=1}^{r}\rho_{k}p_{i}\left(  \omega_{k},\bar{\omega}_{k}\right)
\left|  f_{j}\left(  \omega_{k},\bar{\omega}_{k}\right)  \right|  ^{2}\\
&  =\rho_{j}p_{i}\left(  \omega_{j},\bar{\omega}_{j}\right)  .
\end{align*}
Since $\rho_{j}>0$, then $p_{i}\left(  \omega_{j},\bar{\omega}_{j}\right)
\geq0$. \ Repeating the preceding argument for $1\leq i\leq m$ and $1\leq
j\leq r$, we conclude that $\operatorname*{supp}\nu\subseteq K_{\mathcal{P}}$.

We now count the atoms of $\nu$ that lie in $\mathcal{Z}(p_{i})$. \ Equations
(\ref{eqExi.2}) and (\ref{eq.Mom2}) show that $M_{p_{i}}^{d}\left(
\infty\right)  $ is the moment matrix corresponding to the measure
$p_{i}\,d\nu$, i.e., $M_{p_{i}}^{d}\left(  \infty\right)  =M^{d}\left(
\infty\right)  \left[  p_{i}\,d\nu\right]  $. \ Thus, \cite[Proposition
7.6]{tcmp1} implies that
\[
\operatorname*{card}\operatorname*{supp}(p_{i}d\nu)=\operatorname*{rank}%
M_{p_{i}}^{d}\left(  \infty\right)  =\operatorname*{rank}M_{p_{i}}^{d}\left(
n+k_{i}\right)  .\
\]
We have
\begin{align*}
\Delta_{i}  &  :=\operatorname*{rank}M^{d}\left(  n\right)  (\gamma
)-\operatorname*{rank}M_{p_{i}}^{d}\left(  n+k_{i}\right) \\
&  =\operatorname*{card}\operatorname*{supp}\nu-\operatorname*{card}%
\operatorname*{supp}(p_{i}d\nu)\\
&  =\operatorname*{card}\left(  \operatorname*{supp}\nu\cap\mathcal{Z}\left(
p_{i}\right)  \right)  ,
\end{align*}
whence $\nu$ has precisely $\Delta_{i}$ atoms in $\mathcal{Z}\left(
p_{i}\right)  \;(1\leq i\leq m)$.

For the converse direction, suppose $\nu$ is an $r$-atomic representing
measure for $\gamma$ with $\operatorname*{supp}\nu\subseteq K_{\mathcal{P}}$.
\ Since $\nu$ is a representing measure for $M(\infty)\equiv M^{d}%
(\infty)\left[  \nu\right]  $, \cite[Proposition 7.6]{tcmp1} implies that
\begin{align*}
r  &  =\operatorname{card}\operatorname*{supp}\nu=\operatorname{rank}%
M(\infty)\\
&  \geq\operatorname*{rank}M(n+1)[\nu]\geq\operatorname*{rank}M(n)[\nu
]=\operatorname*{rank}M(n)(\gamma)=r.\
\end{align*}
In particular, $M^{d}\left(  n+1\right)  \left[  \nu\right]  $ is a flat
extension of $M^{d}\left(  n\right)  \left(  \gamma\right)  $, as is
\newline $M^{d}\left(  n+k_{i}\right)  \left[  \nu\right]  \;(1\leq i\leq m)$;
since $\nu$ is also a representing measure for $M^{d}\left(  n+k_{i}\right)
\left[  \nu\right]  $ and $\operatorname*{supp}\nu\subseteq K_{p_{i}}$, then
(\ref{eqLoc.1}) implies that $M_{p_{i}}^{d}\left(  n+k_{i}\right)  \left[
\nu\right]  \geq0\;(1\leq i\leq m)$.
\end{proof}

We next prove Theorem \ref{IntroTheorem1}, the analogue of Theorem
\ref{ThmExi.1} for moment problems on $\mathbb{R}^{N}$. \ We consider a real
$N$-dimensional sequence of degree $2n$, $\beta\equiv\beta^{(2n)}=\{\beta
_{i}\}_{i\in Z_{+}^{N},\left|  i\right|  \leq2n}$, and its moment matrix
$\mathcal{M}\equiv\mathcal{M}^{N}(n)(\beta)$. \ Recall from Theorem
\ref{newthm} that $\beta$ admits a $\operatorname*{rank}\mathcal{M}$-atomic
(minimal) representing measure if and only if $\mathcal{M}\geq0$ and
$\mathcal{M}$ admits a flat moment matrix extension $\mathcal{M}^{N}(n+1)$,
which in turn admits unique successive flat extensions $\mathcal{M}^{N}(n+2)$,
$\mathcal{M}^{N}(n+3)$,... . \ For the reader's convenience, we restate
Theorem \ref{IntroTheorem1}, as follows.

\begin{theorem}
\label{ThmExi.2}Let $\beta\equiv\beta^{(2n)}=\{\beta_{i}\}_{i\in\mathbb{Z}%
_{+}^{N},\left|  i\right|  \leq2n}$ be an $N$-dimensional real sequence, and
let $\mathcal{Q}\equiv\{q_{i}\}_{i=1}^{m}\subseteq\mathbb{C}^{N}[t]$, with
$\deg q_{i}=2k_{i}$ or $2k_{i}-1\;(1\leq i\leq m)$. \ Let $\mathcal{M}%
:=\mathcal{M}^{N}(n)(\beta)$ and let $r:=\operatorname*{rank}\mathcal{M}$.
\ There exists a (minimal) $r$-atomic representing measure for $\mathcal{M}$
supported in $K_{\mathcal{Q}}$ if and only if $\mathcal{M}\geq0$ and
$\mathcal{M}$ admits a flat extension $\mathcal{M}(n+1)$ such that
$\mathcal{M}_{q_{i}}(n+k_{i})\geq0$ $(1\leq i\leq m)$. \ In this case,
$\mathcal{M}(n+1)$ admits a unique representing measure $\mu$, which is an
$r$-atomic (minimal) $K_{\mathcal{Q}}$-representing measure for $\beta$;
moreover, $\mu$ has precisely $r-\operatorname*{rank}\mathcal{M}_{q_{i}%
}(n+k_{i})$ atoms in $\mathcal{Z}(q_{i})\equiv\left\{  t\in\mathbb{R}%
^{N}:q_{i}(t)=0\right\}  $, $1\leq i\leq m$.
\end{theorem}

\begin{proof}
Suppose $\mu$ is a $\operatorname*{rank}\mathcal{M}$-atomic representing
measure for $\beta$ with $\operatorname*{supp}\mu\subseteq K_{\mathcal{Q}}$.
\ Exactly as in the proof of Theorem \ref{newthm} (or of Theorem
\ref{ThmExi.1}), $\mathcal{M}(n+1)[\mu]$ is a flat extension of $\mathcal{M\;}%
(=\mathcal{M}^{N}(n)[\mu]\geq0)$, with unique successive flat extensions
$\mathcal{M}^{N}(n+2)[\mu],\mathcal{M}^{N}(n+3)[\mu],...$ \ \ Since
$\operatorname*{supp}\mu\subseteq K_{\mathcal{Q}}$, for each $f\in
\mathbb{C}_{n}^{N}[t]$, we have $\left\langle \mathcal{M}_{q_{i}}^{N}%
(n+k_{i})[\mu]\hat{f},\hat{f}\right\rangle =\int q_{i}\left|  f\right|
^{2}\;d\mu\geq0$, whence $\mathcal{M}_{q_{i}}^{N}(n+k_{i})[\mu]\geq0\;\;(1\leq
i\leq m)$.

For the converse and the location of the atoms, we first consider the case
when $N$ is even, say $N=2d$. Suppose $\mathcal{M}\equiv\mathcal{M}%
^{2d}(n)(\beta)$ is positive and has a flat extension $\mathcal{M}%
^{2d}(n+1)(\tilde{\beta})$ for which $\mathcal{M}_{q_{i}}^{2d}(n+k_{i}%
)(\tilde{\beta})\geq0\;(1\leq i\leq m)$ (cf. Theorem \ref{newthm}). \ Using
Proposition \ref{prop15} (and as in the proof of the ``even'' case of Theorem
\ref{newthm}), $\mathcal{M}$ corresponds to a complex moment matrix
$M^{d}(n)(\gamma)\;(=L^{(n)\ast}\mathcal{M}L^{(n)})$, and the successive flat
extensions $\mathcal{M}^{2d}(n+j)(\tilde{\beta})$ of $\mathcal{M}$ correspond
to successive flat moment matrix extensions of $M^{d}(n)(\gamma)$ defined by
$M^{d}(n+j)(\tilde{\gamma}):=L^{(n+j)\ast}\mathcal{M}^{2d}(n+j)(\tilde{\beta
})L^{(n+j)}\;(j\geq1)$ (cf. Proposition \ref{prop15}).

We will show that $M_{p_{i}}^{d}\left(  n+k_{i}\right)  (\tilde{\gamma})\geq
0$, where\ $p_{i}:=q_{i}\circ\tau\in\mathbb{C}_{2k_{i}}^{d}[z,\bar
{z}])\;(1\leq i\leq m)$. \ To this end, recall from Lemma \ref{lem10} that
$\tau(z,\bar{z})=(x,y)$, so that $q_{i}=p_{i}\circ\psi\;(1\leq i\leq m)$;
further Proposition \ref{prop15}(vi) implies that
\begin{equation}
\Lambda_{\tilde{\gamma}}(p)=\Lambda_{\tilde{\beta}}(p\circ\psi)\;(p\in
\mathbb{C}^{d}[z,\bar{z}]). \label{eq.Exi.3}%
\end{equation}
We assert that
\begin{equation}
M_{p_{i}}^{d}\left(  n+k_{i}\right)  (\tilde{\gamma})=L^{(n)\ast}%
\mathcal{M}_{q_{i}}^{2d}(n+k_{i})(\tilde{\beta})L^{(n)}\;(1\leq i\leq m).
\label{eq.Exi.4}%
\end{equation}
Indeed, for $f,g\in\mathbb{C}_{n}^{d}[z,\bar{z}])$ and $1\leq i\leq m$, we
have\
\begin{align*}
\left\langle M_{p_{i}}^{d}\left(  n+k_{i}\right)  (\tilde{\gamma})\hat
{f},\hat{g}\right\rangle  &  =\Lambda_{\tilde{\gamma}}(p_{i}f\bar{g})\\
&  =\Lambda_{\tilde{\beta}}((p_{i}f\bar{g})\circ\psi)\;\;\text{(by
(\ref{eq.Exi.3}))}\\
&  =\left\langle \mathcal{M}_{q_{i}}^{2d}(n+k_{i})(\tilde{\beta}%
)\widetilde{f\circ\psi},\widetilde{g\circ\psi}\right\rangle \\
&  =\left\langle \mathcal{M}_{q_{i}}^{2d}(n+k_{i})(\tilde{\beta})L^{(n)}%
\hat{f},L^{(n)}\hat{g}\right\rangle \;\text{(by Lemma \ref{lem10})}\\
&  =\left\langle L^{(n)\ast}\mathcal{M}_{q_{i}}^{2d}(n+k_{i})(\tilde{\beta
})L^{(n)}\hat{f},\hat{g}\right\rangle ,
\end{align*}
whence (\ref{eq.Exi.4}) follows. \ Since $\mathcal{M}_{q_{i}}^{2d}%
(n+k_{i})(\tilde{\beta})\geq0$, then (\ref{eq.Exi.4}) implies that $M_{p_{i}%
}^{d}\left(  n+k_{i}\right)  (\tilde{\gamma})\geq0\;(1\leq i\leq m)$.

Theorem \ref{ThmExi.1} now implies that $\gamma$ has a $\operatorname*{rank}%
M^{d}(n)(\gamma)$-atomic representing measure $\omega$, supported in
$K_{\mathcal{P}}$, and Proposition \ref{PropositionEquiv} shows that $\omega$
corresponds to a $\operatorname*{rank}\mathcal{M}$-atomic representing measure
for $\beta$ supported in $K_{\mathcal{Q}}$. \ Theorem \ref{ThmExi.1} also
implies that $M^{d}\left(  n+1\right)  (\tilde{\gamma})$ admits a unique
representing measure $\nu$, which is a $\operatorname*{rank}M^{d}(n)$-atomic
$K_{\mathcal{P}}$-representing measure for $\gamma$ having
$\operatorname*{rank}M^{d}\left(  n\right)  (\gamma)-\operatorname*{rank}%
M_{p_{i}}^{d}\left(  n+k_{i}\right)  (\tilde{\gamma})$ atoms in $\mathcal{Z}%
(p_{i})\;(1\leq i\leq m)$. \ From Proposition \ref{PropositionEquiv}, $\nu$
corresponds to a unique representing measure $\mu:=\nu\circ\psi$ for
$\mathcal{M}^{2d}(n+1)(\tilde{\beta})$. \ Since, from Proposition
\ref{PropositionEquiv}, $\operatorname*{supp}\nu=\psi(\operatorname*{supp}%
\mu)$, $\mathcal{Z}(p_{i})=\psi\circ\mathcal{Z}(q_{i})$, and
$\operatorname*{rank}M^{d}(n)(\gamma)=\operatorname*{rank}\mathcal{M}$, and
since (\ref{eq.Exi.4}) implies $\operatorname*{rank}M_{p_{i}}^{d}\left(
n+k_{i}\right)  (\tilde{\gamma})=\operatorname*{rank}\mathcal{M}_{q_{i}}%
^{2d}\left(  n+k_{i}\right)  (\tilde{\beta})$, it follows that
$\operatorname*{supp}\mu\subseteq K_{\mathcal{Q}}$ and that $\mu$ has
precisely $\operatorname*{rank}\mathcal{M}-\operatorname*{rank}\mathcal{M}%
_{q_{i}}^{N}\left(  n+k_{i}\right)  (\tilde{\beta})$ atoms in $\mathcal{Z}%
(q_{i})\;(1\leq i\leq m)$. \ The proof of the ``even'' case is now complete.

We now consider the case $N=2d-1$. \ Suppose $\mathcal{M}\equiv\mathcal{M}%
^{2d-1}(n)(\beta)$ is positive and has a flat extension $\mathcal{M}%
^{2d-1}(n+1)(\tilde{\beta})$, with unique successive flat extensions
$\mathcal{M}^{2d-1}(n+j)(\tilde{\beta})\;(j\geq2)$ (cf. Theorem \ref{newthm});
we are assuming $\mathcal{M}_{q_{i}}^{2d-1}(n+k_{i})(\tilde{\beta}%
)\geq0\;(1\leq i\leq m)$. \ As in the proof of the ``odd'' case of Theorem
\ref{newthm}, $\mathcal{M}$ corresponds to the positive moment matrix
$\mathcal{M}^{\smile}\equiv\mathcal{M}^{2d}(n)(\breve{\beta})$, which has a
sequence of successive flat extensions $\mathcal{M}^{2d}(n+j)(\tilde{\lambda
})$ satisfying $\mathcal{M}^{2d}(n+j)(\tilde{\lambda})=\mathcal{\breve{M}%
}^{2d-1}(n+j)(\tilde{\beta})\;(j\geq1)$; the moments of $\tilde{\beta}$ are
related to those of $\tilde{\lambda}$ as in (\ref{betadef}).

Fix $\ell$, $1\leq\ell\leq m$; for $q_{\ell}\equiv\sum b_{\ell,s}t^{s}%
\in\mathbb{C}^{N}[t]$ we let $\breve{q}_{\ell}\in\mathbb{C}^{N+1}[t,u]$ be
given by $\breve{q}_{\ell}(t,u):=q_{\ell}(t)\;(t\in\mathbb{R}^{2d-1}%
,u\in\mathbb{R)}$. \ We claim that $\mathcal{M}_{\breve{q}_{\ell}}%
^{2d}(n+k_{\ell})(\tilde{\lambda})\geq0$. \ To this end, for $i\in
\mathbb{Z}_{+}^{2d-1}$, $j\in\mathbb{Z}_{+}$, recall that $\breve{\imath
}:=(i,j)\in\mathbb{Z}_{+}^{2d}$, and for $t\in\mathbb{R}^{2d-1}$,
$u\in\mathbb{R}$, $\breve{t}:=(t,u)\in\mathbb{R}^{2d}$, so $\breve{t}%
^{\breve{\imath}}=t^{i}u^{j}$. \ We denote $f\in\mathbb{C}_{n}[t,u]$ by
$f(\breve{t})=\sum_{\left|  \breve{\imath}\right|  \leq n}a_{\breve{\imath}%
}\breve{t}^{\breve{\imath}}$, and we define $[f]\in\mathbb{C}_{n}[t]$ by
$[f](t):=\sum_{\left|  \breve{\imath}\right|  \leq n,j=0}a_{\breve{\imath}%
}t^{i}$. \ Now, for $f\in\mathbb{C}_{n}[t,u]$,
\begin{align*}
\left\langle \mathcal{M}_{\breve{q}_{\ell}}^{2d}(n+k_{\ell})(\tilde{\lambda
})\hat{f},\hat{f}\right\rangle  &  =\Lambda_{\tilde{\lambda}}(\breve{q}_{\ell
}\left|  f\right|  ^{2})\\
&  =\sum_{\left|  s\right|  \leq\deg q_{\ell},\left|  \breve{\imath}\right|
,\left|  \breve{\imath}^{\prime}\right|  \leq n}b_{\ell,s}a_{\breve{\imath}%
}\bar{a}_{\breve{\imath}^{\prime}}\tilde{\lambda}_{(s+i+i^{\prime}%
,j+j^{\prime})}\\
&  =\sum_{\left|  s\right|  \leq\deg q_{\ell},\left|  \breve{\imath}\right|
,\left|  \breve{\imath}^{\prime}\right|  \leq n,j=j^{\prime}=0}b_{\ell
,s}a_{\breve{\imath}}\bar{a}_{\breve{\imath}^{\prime}}\tilde{\beta
}_{s+i+i^{\prime}}\;\;\text{(by Remark \ref{Exi.Rem4})}\\
&  =\Lambda_{\tilde{\beta}}(q_{\ell}\left|  [f]\right|  ^{2})=\left\langle
\mathcal{M}_{q_{\ell}}^{2d-1}(n+k_{\ell})(\tilde{\beta})[\hat{f}],[\hat
{f}]\right\rangle \geq0.
\end{align*}
Since $\mathcal{M}_{\breve{q}_{\ell}}^{2d}(n+k_{\ell})(\tilde{\lambda}%
)\geq0\;(1\leq\ell\leq m)$, by the ``even'' case (and its proof, above),
$\mathcal{M}^{2d}(n+1)(\tilde{\lambda})$ admits a unique representing measure
$\tilde{\mu}$, which is a $\operatorname*{rank}\mathcal{M}^{\smile}$-atomic
$K_{\mathcal{Q}^{\smile}}$-representing measure for $\mathcal{M}^{\smile}$
(where $\mathcal{Q}^{\smallsmile}:=\{\breve{q}_{1},...,\breve{q}_{m}\}$), with
precisely $\operatorname*{rank}\mathcal{M}^{\smile}-\operatorname*{rank}%
\mathcal{M}_{\breve{q}_{i}}^{2d}(n+k_{i})(\tilde{\lambda})$ atoms in
$\mathcal{Z}(\breve{q}_{i})\;(1\leq i\leq m)$. \ Write $\tilde{\mu}\equiv
\sum_{s=1}^{r}\rho_{s}\delta_{(t_{s},u_{s})}$; it follows that $\mu
:=\sum_{s=1}^{r}\rho_{s}\delta_{t_{s}}$ is a representing measure for
$\mathcal{M}^{2d-1}(n+1)(\tilde{\beta})$, with $\operatorname*{supp}%
\mu\subseteq K_{\mathcal{Q}}$ and $\operatorname*{card}(\operatorname*{supp}%
\mu\bigcap\mathcal{Z}(q_{i}))=\operatorname*{rank}\mathcal{M}%
-\operatorname*{rank}M_{q_{i}}^{2d}(n+k_{i})(\tilde{\beta})\;(1\leq i\leq m)$.
\ That $\mu$ is the unique representing measure for $M^{2d-1}(n+1)(\tilde
{\beta})$ follows from Theorem \ref{thmA}.
\end{proof}

\begin{proof}
[Proof of Corollary \textup{\ref{cor1.4}}]Suppose $\mathcal{M}(n)$ admits a
positive extension $\mathcal{M}(n+j)$ which in turn has a flat extension
$\mathcal{M}(n+j+1)$ satisfying $\mathcal{M}_{q_{i}}(n+j+k_{i})\geq0\;(1\leq
i\leq m)$. \ We can apply Theorem \ref{IntroTheorem1} to $\mathcal{M}(n+j)$ to
obtain a finitely atomic $K_{\mathcal{Q}}$-representing measure for
$\mathcal{M}(n+j)$, and hence for $\mathcal{M}(n)$. \ For the converse,
suppose $\mathcal{M}(n)$ has a finitely atomic representing measure $\mu$ with
$\operatorname*{supp}\mu\subseteq K_{\mathcal{Q}}$. \ We will estimate the
minimum value of $j$ necessary to obtain a positive extension $\mathcal{M}%
(n+j)$ having a flat extension $\mathcal{M}(n+j+1)$ (with a corresponding
$K_{\mathcal{Q}}$-representing measure). \ Since $\mu$ is finitely atomic, it
has convergent moments of degree $2n+1$. \ Thus, \cite[Theorem 1.4]{tcmp8}
implies that $\mu$ has an \textit{inside} \textit{cubature rule} $\zeta$ of
degree $2n$, with $s:=\operatorname*{card}\operatorname*{supp}\zeta\leq
1+\dim\mathbb{R}_{2n}^{N}[t]=1+\binom{2n+N}{N}$; in particular, $\zeta$ is a
representing measure for $\mathcal{M}(n)$ and $V:=\operatorname*{supp}%
\zeta\subseteq\operatorname*{supp}\mu\;(\subseteq K_{\mathcal{Q}})$. \ Since
$\operatorname*{card}V=s$, Lagrange interpolation implies that every
real-valued function on $V$ agrees on $V$ with a polynomial in $\mathbb{R}%
_{2(s-1)}^{N}[t]$. \ (Indeed, if $V\equiv\{v_{1},...,v_{s}\}$, let $f_{\ell
}(t):=\frac{\prod_{i=1,...,s;i\neq\ell}\left\|  t-v_{i}\right\|  ^{2}}%
{\prod_{i=1,...,s;i\neq\ell}\left\|  t_{j}-v_{i}\right\|  ^{2}}\in
\mathbb{R}_{2(s-1)}^{N}[t]\;(1\leq\ell\leq s)$. \ Then any function
$f:V\rightarrow\mathbb{R}$ satisfies $f=\sum_{\ell=1}^{s}f(v_{\ell})f_{\ell}%
$.) \ In particular, if $i\in\mathbb{Z}_{+}^{N}$ with $\left|  i\right|
=2s-1$, there exists $p_{i}\in\mathbb{R}_{2(s-1)}^{N}[t]$ such that
$t^{i}-p_{i}(t)|_{V}\equiv0$. \ By Proposition \ref{MomProp1}, $T^{i}%
=p_{i}(T)$ in $\mathcal{C}_{\mathcal{M}(2s-1)[\zeta]}$, and since $\deg
p_{i}<\left|  i\right|  $, it follows that $\mathcal{M}(2s-1)[\zeta]$ is a
flat extension of $\mathcal{M}(2s-2)[\zeta]$. \ Theorem \ref{newthm} implies
that $\mathcal{M}(2s-1)[\zeta]$ has unique successive flat moment matrix
extensions, and it is clear from the preceding argument that these extensions
are $\mathcal{M}(2s)[\zeta],\mathcal{M}(2s+1)[\zeta],...$ . \ Since
$V\subseteq K_{\mathcal{Q}}$, it follows immediately that $\mathcal{M}_{q_{i}%
}(2s-2+k_{i})[\zeta]\geq0\;\;(1\leq i\leq m)$. \ If $n\leq2(s-1)$, then
$j:=2(s-1)-n$ satisfies our requirements, and $j\leq2\binom{2n+N}{N}-n$. \ If
$n>2(s-1)$, then $\mathcal{M}(n)=\mathcal{M}(2s-1)[\zeta]$ or $\mathcal{M}(n)$
is one of the successive extensions of $\mathcal{M}(2s-1)[\zeta]$ listed
above, and in this case we can take $j:=0$.
\end{proof}

\begin{remark}
(i) In the case $N=2,n>2$, the estimate for $j$ ($j\leq4n^{2}+5n+2$) can be
improved to $j\leq2n^{2}+6n+6$ (cf. \cite[Theorem 1.5]{tcmp3}). \ We also note
that in several examples that we have studied which require $j>0$, the flat
extension $\mathcal{M}(n+j+1)$ can be realized with $j=1$ \cite{tcmp7},
\cite{tcmp9}, \cite{FiaOT}, \cite{FiPe}. \ In particular, if $K_{\mathcal{Q}}$
is a degenerate hyperbola and $\mathcal{M}(n)$ has a $K_{\mathcal{Q}}%
$-representing measure, $\mathcal{M}(n)$ might not have a flat extension, but
in this case there is always a positive extension $\mathcal{M}(n+1)$ that has
a flat extension $\mathcal{M}(n+2)$ \cite{tcmp9}.\newline (ii) Corollary
\ref{cor1.4} implies an exact analogue for complex moment sequences.
\end{remark}

\begin{proof}
[Added in Proof]After completing this paper we learned of recent related
papers of Professor Monique Laurent (\cite{Lau1}, \cite{Lau2}, \cite{Lau3}).
\ In \cite[Theorem 1.2]{Lau3}, Laurent gives an alternate proof of Corollary
\ref{newcor}, using algebraic techniques (e.g., Nullstellensatz) to prove the
existence of a unique, $r$-atomic, representing measure corresponding to a
rank $r$ positive infinite moment matrix. \ Using this result and Theorem
\ref{newthm}, Laurent then provides a short proof of Theorem
\ref{IntroTheorem1} (\cite[Theorem 1.6]{Lau3}). \ This proof is based on
general properties of localizing matrices, and circumvents our explicit
calculations of localizing matrices in Section 4. \ For applications, to
explicitly compute representing measures, it is still necessary to use the
computational formula of Theorem \ref{ThmLoc.2}.
\end{proof}

\end{document}